\documentclass[10.5pt]{article}

\usepackage{amsmath}
\usepackage{amscd}
\usepackage{amssymb}
\usepackage{rotating}
\usepackage{amsfonts}
\usepackage{amsthm}
\usepackage{verbatim}
\usepackage{stmaryrd}
\usepackage[top=1in, bottom=1in, left=1.25in, right=1.25in]{geometry}

\usepackage{titlesec}
       \titleformat{\chapter}[display]
             {\normalfont\Large\bfseries}{\thechapter}{11pt}{\Large}
       \titleformat{\section}
             {\normalfont\large\bfseries}{\thesection}{10.5pt}{\large}
       \titlespacing*{\chapter}{0pt}{0pt}{15pt} 
       \titlespacing*{\section}{0pt}{3.5ex plus 1ex minus .2ex}{2.3ex plus .2ex}

    \usepackage[titletoc]{appendix}

\input xy
\xyoption{all}

\newcommand{\cM}{{\mathcal M}}

\newcommand{\Mbar}{\overline{\cM}}

\newcommand{\pqed}{\hfill\qedsymbol\\}
\newcommand{\nn}{\nonumber}
\newcommand{\Res}{\mathrm{Res}}

\newtheorem{theorem}{Theorem}[section]
\newtheorem{theorem/definition}{Theorem/Definition}[section]
\newtheorem{Theorem}{Theorem}
\newtheorem{proposition}{Proposition}[section]
\newtheorem{lemma}{Lemma}[section]

\newtheorem{Conjecture}{Conjecture}

\theoremstyle{remark}
\newtheorem{remark}{Remark}[section]

\theoremstyle{definition}

\newcommand{\be}{\begin{equation}}
\newcommand{\ee}{\end{equation}}
\newcommand{\bea}{\begin{eqnarray}}
\newcommand{\eea}{\end{eqnarray}}
\newcommand{\ben}{\begin{eqnarray*}}
\newcommand{\een}{\end{eqnarray*}}
\newcommand{\bet}{\begin{equation}
\begin{split}}
\newcommand{\eet}{\end{split}
\end{equation}}

\begin{document}
\title
{A Conjectural Formula for Genus One Gromov-Witten Invariants of a Class of Local Calabi-Yau $n$-folds}
\author{\normalsize Xiaowen Hu}
\maketitle

\begin{abstract}
We conjecture a formula for the generating function of genus one Gromov-Witten invariants of the local Calabi-Yau manifolds which are the total spaces of splitting bundles over projective spaces. We prove this conjecture in several special cases, and assuming the validity of our conjecture we check the integrality of genus one BPS numbers of local Calabi-Yau 5-folds defined by A. Klemm and R. Pandharipande.
\end{abstract}
\section{Introduction}
After a series of splendid works with Jun Li and R.Vakil (see \cite{Zinger1} the references therein), A.Zinger finally explicitly computed the genus one Gromov-Witten invariants of Calabi-Yau hypersurfaces in projective spaces. This result is generalized to complete intersections in projective spaces by A.Popa in (\cite{Popa1}). Our object is to find a similar formula for the genus one Gromov-Witten invariants of the local Calabi-Yau $n$-fold
\bea\label{56}
X=\mathrm{Tot}\big(\mathcal{O}(-c_{1})\oplus\cdots \mathcal{O}(-c_{m})\rightarrow \mathbb{P}^{n-m}\big),
\eea
 where $c_{i}\in \mathbb{Z}^{>0}$ for $1\leq i\leq m$ and $\sum_{i=1}^{m}c_{i}=n-m+1$.\\
  Let us first recall Zinger's formula. Let the target space $Y$ be a degree $n$ hypersurface in $\mathbb{P}^{n-1}$. For $q=0,1,...$, define $I_{0,q}$ by
\bea\label{51}
\sum_{q=0}^{\infty}I_{0,q}(t)w^{q}=e^{wt}\sum_{d=0}^{\infty}e^{dt}\frac{\prod_{r=1}^{nd}(nw+r)}{\prod_{r=1}^{d}(w+r)^{n}}.
\eea
It is easy to see that for $0\leq q\leq n-2$, $I_{0,q}$ are solutions of the Picard-Fuchs operator
\bea\label{52}
\mathcal{L}=\Big(\frac{d}{dt}\Big)^{n-1}-ne^{t}\prod_{r=1}^{n-1}(n\frac{d}{dt}+r).
\eea
For $q\geq p\geq 0$, we inductively define
\bea\label{53}
I_{p,q}(t)=\frac{d}{dt}\Big(\frac{I_{p-1,q}(t)}{I_{p-1,p-1}(t)}\Big),
\eea
and the mirror map is given by
\bea\label{54}
T=\frac{I_{0,1}(t)}{I_{0,0}(t)}.
\eea
Thus $T-t$ and $I_{p,p}(t)$  are series of $e^{t}$ for $p\geq 0$.
The genus one degree $d$ Gromov-Witten invariants $N_{1,d}^{X}$ are given by
\bea\label{55}
\sum_{d=1}^{\infty}N_{1,d}^{Y}e^{dT}&=&\Big(\frac{(n-2)(n+1)}{48}+\frac{1-(1-n)^{n}}{24n^{2}}\Big)(T-t)+\frac{n^{2}-1+(1-n)^{n}}{24n}\ln I_{0,0}(t)\nn\\
&&-\left\{\begin{array}{lll}
\frac{n-1}{48}\ln(1-n^{n}e^{t})+\sum_{p=0}^{\frac{n-3}{2}}\frac{(n-1-2p)^{2}}{8}\ln I_{p,p}(t), & if & 2\nmid n;\\
\frac{n-4}{48}\ln(1-n^{n}e^{t})+\sum_{p=0}^{\frac{n-4}{2}}\frac{(n-2p)(n-2-2p)}{8}\ln I_{p,p}(t), & if & 2\mid n.\\
\end{array}
\right.\nn\\
\eea
Before Zinger's work, the formula for $n=5$ or 6 ($Y$ is a quintic 3-fold or a sextic 4-fold, resp.) had been conjectured via mirror symmetry and physical arguments on the B-side, see \cite{BCOV1} and \cite{KP}. For $n\geq 7$, the B-side interpretation is still absent, at least to the best knowledge of the author. \\

The Gromov-Witten invariants of local Calabi-Yau manifolds which are total spaces of vector bundles over toric varieties are in principle less difficult to compute, because we can directly apply the virtual localization method. But in dimension greater than 3, it seems not easy to get a closed formula due to the complicated combinatorics. Thus to get a formula for local Calabi-Yau spaces, a possible approach is just to adapt Zinger's method to the local case, i.e., we need to\\
1. Find a \emph{standard vs reduced} comparison formula for relevant Hodge integrals on $\Mbar_{1,k}(\mathbb{P}^{n-m+1},d)$ and $\overline{\mathfrak{M}}_{1,k}^{0}(\mathbb{P}^{n-m+1},d)$.\\
2. Find a formula for Hodge integrals on $\widetilde{\mathcal{M}}_{1,k}$.\\
3. Write the Hodge integrals on $\overline{\mathfrak{M}}_{1,k}^{0}(\mathbb{P}^{n-m+1},d)$ as contributions of graphs by localization.\\
4. Generalize the combinatorial arguments in \cite{Zinger1} to the local cases.\\

In principle also, the above procedure should be less difficult than that of the compact cases, since in the latter cases the involved sheaves 
$R^{0}\pi_{*}f^{*}\mathcal{O}(n-m+2)$ is not locally free. We have made some progress on this and hope to address it in the future. In this article, however, we get a formula by a mixture of physical arguments and mathematical observations on Zinger's proof, and we check the formula by proving it in several most simple cases, and also by checking the integrality of the BPS numbers of local Calabi-Yau 5-folds.\\

Now let us take a closer look at (\ref{55}). For the first term, the coefficient of $T-t$ physically (see \cite{BCOV1}) comes from the integral
\bea
\frac{1}{24}\int_{Y}k\wedge c_{n-3}(Y),
\eea
where\footnote{The potential $\mathcal{F}_{1}$ differs from the nowadays usual choice of potential by a factor 2, so the coefficient $\frac{1}{12}$ is taken as $\frac{1}{24}$ here. } $k$ is the  K\"{a}hler class of $Y$ associated with the variable $T$, and  is $H$ here,  the class induced by the hyperplane class in the ambient space $\mathbb{P}^{n-1}$. The Chern class is easily computed
\bea
c_{n-3}(Y)=\Bigg(\frac{(n-2)(n+1)}{2n}+\frac{1-(1-n)^{n}}{n^{3}}\Bigg)H^{n-3}.
\eea
For the local case, for the target space $X$ of the form (\ref{56}), the series corresponding to (\ref{51}) is
\bea\label{57}
\sum_{q=0}^{\infty}I_{0,q}(t)w^{q}=e^{wt}\sum_{d=0}^{\infty}e^{dt}\frac{\prod_{i=1}^{m}\prod_{s=0}^{c_{i}d-1}(-c_{i}w-s)}
{\prod_{s=1}^{d}(w+s)^{n-m+1}},
\eea
which encodes the genus zero one-point and two-point Gromov-Witten invariants of $X$ by \cite{Popa2}.
It is easy to see that, when $m>1$ the mirror map is the identity map $T=t$, so the first term of (\ref{55}) has no counterpart in these cases. When $m=1$, $X$ is the total space of the canonical bundle of $\mathbb{P}^{n-1}$, and $c_{n-1}(X)=-\frac{n(n+1)(n-2)}{2}H^{n-1}$. The K\"{a}hlar class is still $H$, but the integral of $c_{n-1}(X)\wedge H$ over the local space $X$ should be taken as the integral of the (formal) quotient of $c_{n-1}(X)\wedge H$ by the  Euler class of $\mathcal{O}(-n)$ over the compact part $\mathbb{P}^{n-1}$, as a general principle\footnote{Writing the local Gromov-Witten invariants as Hodge integrals over the moduli space of stable maps to the compact part, to make the WDVV equation still hold, we need to cancel one of the two copies of contributions of the Euler class at the node, in the usual derivation of the WDVV equation. }.\\
We can also  get the same result in another way. In the mathematical proof of Zinger, the  coefficient of $T-t$ comes from a computation of residues. In fact, the first term of the coefficient comes from a residue at 0, and the second term from a residue at $-n$. In the local case, by a speculation on Zinger's proof, there should be no residues at $-n$ and the residue at $0$ is the same as the global case. So the counterpart of the first term in the formula for $K_{\mathbb{P}^{n-1}}$ should be
\bea
\frac{(n+1)(n-2)}{48}(T-t).
\eea
For the second term of (\ref{55}), since in the local case we always have $I_{0,0}(t)=1$ from (\ref{57}), it has no counterpart in the local case.\\

For the third term of (\ref{55}), we follow the arguments in \cite{KP}. By some physical argument, this term comes from the behavior of the potential at the \emph{conifold} point of the moduli space on the B-side, and
the coefficient $-\frac{n-1}{48}$ or $-\frac{n-4}{48}$ ($n$ is odd or even, resp.) should be universal. The  $1-n^{n}e^{t}$ comes from the discriminant of the Picard-Fuchs operator (\ref{52}). In the local case, the Picard-Fuchs operator is
\ben
\mathcal{L}=\Big(\frac{d}{dt}\Big)^{n-m+1}-e^{t}\prod_{i=1}^{m}\prod_{s=0}^{c_{i}-1}(-c_{i}\frac{d}{dt}-s),
\een
and the discriminant is
\bea
\Delta=1-\prod_{i=1}^{m}(-c_{i})^{c_{i}}e^{t}.
\eea
So the counterpart of the third term in the local case should be
\bea
-\left\{\begin{array}{lll}
\frac{n+1}{48}\ln(1-\prod_{i=1}^{m}(-c_{i})^{c_{i}}e^{t}),& \mathrm{if} &2\nmid n;\\
\frac{n-2}{48}\ln(1-\prod_{i=1}^{m}(-c_{i})^{c_{i}}e^{t}),& \mathrm{if}& 2\mid n.
\end{array}\right.
\eea
The fourth group of terms of (\ref{55}) seems the most mysterious.  On one hand, I believe that, to get a series of $e^{t}$ (not a mixture of $t$ and $e^{t}$, or equivalently, without $\log q$ terms, where $q=e^{t}$) from the solutions of the corresponding Picard-Fuchs equation, and to encode enough data from these solutions to get the genus one invariants, the inductive procedure (\ref{53}) is somewhat \emph{ubiquitous}, and thus in the same way we obtain $I_{p,p}(t)$ in the local case.
On the other hand, by a speculation on the argument in \cite{KP}, I believe that if one could find a B-side interpretation of (\ref{55}), the coefficient of $I_{p,p}(t)$ would come from the fact $h^{p,p}=1$ for $0\leq p\leq n-2$ (corresponding to the Ramond-Ramond sector on the B-side) and the elementary identities
\ben
\frac{1}{2}+\frac{3}{2}+\cdots+\frac{n-2-2p}{2}=\frac{(n-1-2p)^{2}}{8}
\een
or
\ben
\frac{2}{2}+\frac{4}{2}+\cdots+\frac{n-2-2p}{2}=\frac{(n-2p)(n-2-2p)}{8}
\een
for $n$ is odd or even, resp.. So the counterpart of the fourth group of  terms in the local cases should be
\bea
-\left\{\begin{array}{lll}
\sum_{p=1}^{(n-1)/2}\frac{(n+1-2p)^{2}}{8}\ln I_{p,p}(t),& \mathrm{if} &2\nmid n;\\
\sum_{p=1}^{(n-2)/2}\frac{(n+2-2p)(n-2p)}{8}\ln I_{p,p}(t),& \mathrm{if}& 2\mid n.
\end{array}\right.
\eea
Combining the above discussions, we obtain the following
\begin{Conjecture}\label{60}
Let $X$ be of the form (\ref{56}). For $m=1$, we have
\bea\label{58}
\sum_{d=1}^{\infty}N_{1,d}^{X}e^{dT}=\frac{(n-2)(n+1)}{48}(T-t)-\left\{\begin{array}{ll}
\frac{n+1}{48}\ln(1+n^{n}e^{t})+\sum_{p=1}^{(n-1)/2}\frac{(n+1-2p)^{2}}{8}\ln I_{p,p}(t),& \mathrm{if}\hspace{0.2cm}2\nmid n;\\
\frac{n-2}{48}\ln(1-n^{n}e^{t})+\sum_{p=1}^{(n-2)/2}\frac{(n+2-2p)(n-2p)}{8}\ln I_{p,p}(t),& \mathrm{if}\hspace{0.2cm}2\mid n.
\end{array}\right.\nn\\
\eea
For $m\geq 2$, we have (in these cases $T=t$)
\bea\label{59}
\sum_{d=1}^{\infty}N_{1,d}^{X}e^{dt}=-\left\{\begin{array}{lll}
\frac{n+1}{48}\ln(1-\prod_{i=1}^{m}(-c_{i})^{c_{i}}e^{t})+\sum_{p=1}^{(n-1)/2}\frac{(n+1-2p)^{2}}{8}\ln I_{p,p}(t),& \mathrm{if} &2\nmid n;\\
\frac{n-2}{48}\ln(1-\prod_{i=1}^{m}(-c_{i})^{c_{i}}e^{t})+\sum_{p=1}^{(n-2)/2}\frac{(n+2-2p)(n-2p)}{8}\ln I_{p,p}(t),& \mathrm{if}& 2\mid n.
\end{array}\right.
\eea
\end{Conjecture}
In fact, the above discussions suggest a \emph{recipe} to get genus one Gromov-Witten invariants from genus zero invariants for Calabi-Yau $n$-folds with $h^{1,1}=1$. Thus one can try to make similar conjectures for, e.g., Calabi-Yau complete intersections in Grassmannians. It is very desirable to give a B-side interpretation of these formulae, e.g., by solving the $tt^{*}$-equations.\\

The $n=3$ and $n=4$ cases of the conjecture \ref{60} has been given in \cite{ABK} and \cite{KP}. The main theorem of this article is
\begin{Theorem}\label{61}
The conjecture \ref{60} holds for degree one invariants, and holds for $X=\mathrm{Tot}\big(\mathcal{O}(-1)^{\oplus (l+1)}\rightarrow \mathbb{P}^{l}\big)$ and $X=\mathrm{Tot}\big(\mathcal{O}(-1)^{\oplus(l-1)}\oplus\mathcal{O}(-2)\rightarrow \mathbb{P}^{l}\big)$ in all degrees, for $l\geq 1$.
\end{Theorem}
We prove this theorem by virtual localization (\cite{GP}). Finally, we check the integrality of $n_{1,d}$ defined for Calabi-Yau 5-folds in \cite{PZ}, from our conjectural formulae (\ref{58}) and (\ref{59}). \\

\noindent\textbf{Conventions}: \\
$\bullet$ We use $[x^{k}]\Big(f(x)\Big)$ to represent the coefficient of $x^{k}$ in the Laurent expansion of $f(x)$ at $x=0$. In this article $x$ may be $q$, $e^{t}$, $Q$ or $w$.\\
$\bullet$ Since the compact part of the target spaces that we consider in this article are always projective spaces, we use $H$ to denote the hyperplane class throughout. Also, $N_{1,d}^{X}$ always denotes the genus one Gromov-Witten invariants of the Calabi-Yau space $X$ with no insertion.\\
$\bullet$ We always understand $Q=e^{T}$ and $q=e^{t}$. In the first three sections we usually use $e^{t}$ and $e^{T}$. In the  section 4 we use $Q$ and $q$, and understand that $I_{p,p}(q)$ means replacing $e^{t}$ by $q$ in the expansion of $I_{p,p}(t)$. \\
$\bullet$ In the graphs that represent the fixed loci in the moduli spaces of genus one stable maps, $\circ$ represents a genus one component, and $\bullet$ represents a genus zero component.\\
$\bullet$ The formal integrals over $\Mbar_{0,1}$ and $\Mbar_{0,2}$ are understood as extending the range of $n$ in the following identity to $n\geq 1$:
\ben
\int_{\Mbar_{0,n}}\frac{1}{\prod_{i=1}^{n}(w_{i}-\psi_{i})}=\frac{1}{\prod_{i=1}^{n}w_{i}}\big(\sum_{i=1}^{n}\frac{1}{w_{i}}\big)^{n-3}.
\een
\textit{Acknowledgements.} The author thanks Prof.~Jian Zhou for his great patience and guidance during all the time. He also thanks Huazhong Ke, Jie Zhou, Xiaobo Zhuang, and Di Yang for helpful discussions. He especially thanks Jie Zhou for carefully reading  an earlier version of the introduction and giving suggestions.

\section{Degree one invariants}
The genus one degree one invariants of local Calabi-Yau $n$-folds of the form of (\ref{56}) can be easily computed by virtual localization.
Let the torus $(\mathbb{C}^{*})^{n-m+1}$ acts on $\mathbb{P}^{n-m}$ with fixed point $P_{i}$, $1\leq i\leq n-m+1$, such that the $n-m$ weights at $P_{i}$  is $\alpha_{i}-\alpha_{k}$, for $k\in\{1,\cdots,n-m+1\}\backslash\{i\}$. We choose the linearizations of $\mathcal{O}(-c_{i})$
with weight $-c_{i}\alpha_{k}$ at $P_{k}$, for $1\leq i\leq m$, $1\leq k\leq n-m+1$. The torus action naturally induces an action on $\Mbar_{1,0}(\mathbb{P}^{n-m+1},1)$, whose fixed loci are corresponding to the graphs of the form

$$ \Gamma_{ij}=\xy
(0,0); (10,0), **@{-};
(0,0)*+{\circ};(10,0)*+{\bullet};(5,2);
(0,3)*+{i};(10,3)*+{j};
\endxy,
$$
where $1\leq i\neq j\leq n-m+1$. Let us first assume $m=1$. Then the contribution of $\Gamma_{ij}$ is
\ben
&&\int_{\Mbar_{1,1}}\frac{(\alpha_{j}-\alpha_{i})\prod_{k\neq i}\Lambda_{1}^{\vee}(\alpha_{i}-\alpha_{k})\cdot \Lambda_{1}^{\vee}(-n\alpha_{i}) \prod_{a=1}^{n-1}\big(-n\alpha_{j}+a(\alpha_{j}-\alpha_{i})\big)}{(\alpha_{i}-\alpha_{j}-\psi)
(\alpha_{i}-\alpha_{j})(\alpha_{j}-\alpha_{i})\prod_{k\neq i,j}\prod_{a=0}^{1}\big(\alpha_{i}-\alpha_{k}+a(\alpha_{j}-\alpha_{i})\big)}\\
&=&\frac{(-1)^{n}n}{24}\frac{\alpha_{i} \prod_{a=1}^{n-1}\big((n-a)\alpha_{j}+a\alpha_{i}\big)}{
\prod_{k\neq j}(\alpha_{j}-\alpha_{k})}
\Big(\sum_{k\neq i,j}\frac{1}{\alpha_{i}-\alpha_{k}}-\frac{1}{n\alpha_{i}}\Big).
\een

Note that
\ben
\sum_{j\neq i}\frac{\prod_{a=1}^{n-1}\big((n-a)\alpha_{j}+a\alpha_{i}\big)}{\prod_{k\neq j}(\alpha_{j}-\alpha_{k})}&=&
n^{n-1}\alpha_{i}^{n-1}\sum_{j\neq i}\frac{1}{(\alpha_{j}-\alpha_{i})\prod_{k\neq i,j}(\alpha_{j}-\alpha_{k})}
+(n-1)!\\
&=&-\frac{n^{n-1}\alpha_{i}^{n-1}}{\prod_{j\neq i}(\alpha_{i}-\alpha_{j})}+(n-1)!,
\een

\ben
\sum_{j\neq i}\frac{\prod_{a=1}^{n-1}\big((n-a)\alpha_{j}+a\alpha_{i}\big)}{(\alpha_{j}-\alpha_{i})^{2}\prod_{k\neq i, j}(\alpha_{j}-\alpha_{k})}&=&
\sum_{j\neq i}\frac{n^{n-1}\alpha_{i}^{n-1}}{(\alpha_{j}-\alpha_{i})^{2}\prod_{k\neq i, j}(\alpha_{j}-\alpha_{k})}
+\sum_{j\neq i}\frac{n^{n-2}\alpha_{i}^{n-2}\cdot\frac{n(n-1)}{2}}{(\alpha_{j}-\alpha_{i})\prod_{k\neq i, j}(\alpha_{j}-\alpha_{k})}\\
&=&\frac{n^{n-1}\alpha_{i}^{n-1}}{\prod_{j\neq i}(\alpha_{i}-\alpha_{j})}\sum_{j\neq i}\frac{1}{\alpha_{i}-\alpha_{j}}-\frac{n^{n-1}(n-1)\alpha_{i}^{n-2}}{2\prod_{j\neq i}(\alpha_{i}-\alpha_{j})},
\een
which are easily to show by the residue theorem on $\mathbb{P}^{1}$.
Thus we have
\ben
&&\sum_{j\neq i}\frac{n\alpha_{i} \prod_{a=1}^{n-1}\big((n-a)\alpha_{j}+a\alpha_{i}\big)}{
\prod_{k\neq j}(\alpha_{j}-\alpha_{k})}
\Big(\sum_{k\neq i,j}\frac{1}{\alpha_{i}-\alpha_{k}}-\frac{1}{n\alpha_{i}}\Big)\\
&=&n\alpha_{i}\Big(\sum_{k\neq i}\frac{1}{\alpha_{i}-\alpha_{k}}\Big)\sum_{j\neq i}\frac{\prod_{a=1}^{n-1}\big((n-a)\alpha_{j}+a\alpha_{i}\big)}{\prod_{k\neq j}(\alpha_{j}-\alpha_{k})}
-\sum_{j\neq i}\frac{\prod_{a=1}^{n-1}\big((n-a)\alpha_{j}+a\alpha_{i}\big)}{\prod_{k\neq j}(\alpha_{j}-\alpha_{k})}\\
&&+n\alpha_{i}\sum_{j\neq i}\frac{\prod_{a=1}^{n-1}\big((n-a)\alpha_{j}+a\alpha_{i}\big)}{(\alpha_{j}-\alpha_{i})^{2}\prod_{k\neq i, j}(\alpha_{j}-\alpha_{k})}\\
&=&\Big[n\alpha_{i}\Big(\sum_{k\neq i}\frac{1}{\alpha_{i}-\alpha_{k}}\Big)-1\Big]\Big(-\frac{n^{n-1}\alpha_{i}^{n-1}}{\prod_{j\neq i}(\alpha_{i}-\alpha_{j})}+(n-1)!\Big)\\
&&+\frac{n^{n}\alpha_{i}^{n}}{\prod_{j\neq i}(\alpha_{i}-\alpha_{j})}\sum_{j\neq i}\frac{1}{\alpha_{i}-\alpha_{j}}
-\frac{n^{n}(n-1)\alpha_{i}^{n-1}}{2\prod_{j\neq i}(\alpha_{i}-\alpha_{j})}\\
&=&n!\alpha_{i}\Big(\sum_{k\neq i}\frac{1}{\alpha_{i}-\alpha_{k}}\Big)-(n-1)!
-\frac{n^{n-1}(n-2)(n+1)\alpha_{i}^{n-1}}{2\prod_{j\neq i}(\alpha_{i}-\alpha_{j})},
\een
and thus
\bea\label{62}
&&\sum_{i=1}^{n}\sum_{j\neq i}\frac{n\alpha_{i} \prod_{a=1}^{n-1}\big((n-a)\alpha_{j}+a\alpha_{i}\big)}{
\prod_{k\neq j}(\alpha_{j}-\alpha_{k})}
\Big(\sum_{k\neq i,j}\frac{1}{\alpha_{i}-\alpha_{k}}-\frac{1}{n\alpha_{i}}\Big)\nn\\
&=&n!\cdot\frac{n(n-1)}{2}-n!
-\frac{n^{n-1}(n-2)(n+1)}{2}=\frac{(n!-n^{n-1})(n-2)(n+1)}{2}.
\eea

Now assume $m\geq 2$.  The contribution of $\Gamma_{ij}$ is
\ben
&&\int_{\Mbar_{1,1}}\frac{(\alpha_{j}-\alpha_{i})\prod_{k\neq i}\Lambda_{1}^{\vee}(\alpha_{i}-\alpha_{k})\cdot \prod_{l=1}^{m}\Big(\Lambda_{1}^{\vee}(-c_{l}\alpha_{i}) \prod_{a=1}^{c_{l}-1}(-c_{l}\alpha_{j}+a(\alpha_{j}-\alpha_{i}))\Big)}
{(\alpha_{i}-\alpha_{j}-\psi)
(\alpha_{i}-\alpha_{j})(\alpha_{j}-\alpha_{i})\prod_{k\neq i,j}(\alpha_{i}-\alpha_{k})(\alpha_{j}-\alpha_{k})}\\
&=&\frac{(-1)^{n-m+1}}{24}\frac{\alpha_{i}^{m}\prod_{l=1}^{m}\Big(c_{l} \prod_{a=1}^{c_{l}-1}\big((c_{l}-a)\alpha_{j}+a\alpha_{i}\big)\Big)}{(\alpha_{j}-\alpha_{i})\prod_{k\neq i,j}(\alpha_{j}-\alpha_{k})}
\Big(\sum_{k\neq i,j}\frac{1}{\alpha_{i}-\alpha_{k}}-\sum_{l=1}^{m}\frac{1}{c_{l}\alpha_{i}}\Big).
\een
Similar to the $m=1$ case, we have
\ben
\sum_{j\neq i}\frac{\prod_{l=1}^{m} \prod_{a=1}^{c_{l}-1}\big((c_{l}-a)\alpha_{j}+a\alpha_{i}\big)}{\prod_{k\neq j}(\alpha_{j}-\alpha_{k})}=-\frac{\prod_{l=1}^{m}c_{l}^{c_{l}-1}\alpha_{i}^{n-2m+1}}{\prod_{j\neq i}(\alpha_{i}-\alpha_{j})},
\een
and
\ben
&&\sum_{j\neq i}\frac{\prod_{l=1}^{m} \prod_{a=1}^{c_{l}-1}\big((c_{l}-a)\alpha_{j}+a\alpha_{i}\big)}
{(\alpha_{j}-\alpha_{i})^{2}\prod_{k\neq i, j}(\alpha_{j}-\alpha_{k})}\\
&=&\frac{\prod_{l=1}^{m}c_{l}^{c_{l}-1}\alpha_{i}^{n-2m+1}}{\prod_{j\neq i}(\alpha_{i}-\alpha_{j})}\sum_{j\neq i}\frac{1}{\alpha_{i}-\alpha_{j}}-\frac{(n-2m+1)\prod_{l=1}^{m}c_{l}^{c_{l}-1}\alpha_{i}^{n-2m}}{2\prod_{j\neq i}(\alpha_{i}-\alpha_{j})}.
\een
So
\ben
&&\sum_{j\neq i}\frac{\alpha_{i}^{m}\prod_{l=1}^{m}\Big(c_{l} \prod_{a=1}^{c_{l}-1}\big((c_{l}-a)\alpha_{j}+a\alpha_{i}\big)\Big)}{(\alpha_{j}-\alpha_{i})\prod_{k\neq i,j}(\alpha_{j}-\alpha_{k})}
\Big(\sum_{k\neq i,j}\frac{1}{\alpha_{i}-\alpha_{k}}-\sum_{l=1}^{m}\frac{1}{c_{l}\alpha_{i}}\Big)\\
&=&\sum_{j\neq i}\Big[\sum_{k\neq i}\frac{1}{\alpha_{i}-\alpha_{k}}\cdot\frac{\alpha_{i}^{m}\prod_{l=1}^{m}\Big(c_{l} \prod_{a=1}^{c_{l}-1}\big((c_{l}-a)\alpha_{j}+a\alpha_{i}\big)\Big)}{(\alpha_{j}-\alpha_{i})\prod_{k\neq i,j}(\alpha_{j}-\alpha_{k})}
+\frac{\alpha_{i}^{m}\prod_{l=1}^{m}\Big(c_{l} \prod_{a=1}^{c_{l}-1}\big((c_{l}-a)\alpha_{j}+a\alpha_{i}\big)\Big)}{(\alpha_{j}-\alpha_{i})^{2}\prod_{k\neq i,j}(\alpha_{j}-\alpha_{k})}\\
&&-\sum_{l=1}^{m}\frac{1}{c_{l}\alpha_{i}}\frac{\alpha_{i}^{m}\prod_{l=1}^{m}\Big(c_{l} \prod_{a=1}^{c_{l}-1}\big((c_{l}-a)\alpha_{j}+a\alpha_{i}\big)\Big)}{(\alpha_{j}-\alpha_{i})\prod_{k\neq i,j}(\alpha_{j}-\alpha_{k})}\Big]\\
&=&-\sum_{k\neq i}\frac{1}{\alpha_{i}-\alpha_{k}}\cdot\frac{\prod_{l=1}^{m}c_{l}^{c_{l}}\alpha_{i}^{n-m+1}}{\prod_{j\neq i}(\alpha_{i}-\alpha_{j})}+\frac{\prod_{l=1}^{m}c_{l}^{c_{l}}\alpha_{i}^{n-m+1}}{\prod_{j\neq i}(\alpha_{i}-\alpha_{j})}\sum_{j\neq i}\frac{1}{\alpha_{i}-\alpha_{j}}\\
&&-\frac{(n-2m+1)\prod_{l=1}^{m}c_{l}^{c_{l}}\alpha_{i}^{n-m}}{2\prod_{j\neq i}(\alpha_{i}-\alpha_{j})}
+\sum_{l=1}^{m}\frac{1}{c_{l}\alpha_{i}}\frac{\prod_{l=1}^{m}c_{l}^{c_{l}}\alpha_{i}^{n-m+1}}{\prod_{j\neq i}(\alpha_{i}-\alpha_{j})}\\
&=&\Big(\sum_{l=1}^{m}\frac{1}{c_{l}}-\frac{n-2m+1}{2}\Big)\frac{\prod_{l=1}^{m}c_{l}^{c_{l}}\alpha_{i}^{n-m}}{\prod_{j\neq i}(\alpha_{i}-\alpha_{j})}.
\een
Therefore for $m\geq 2$ we obtain
\bea\label{63}
N_{1,d}^{X}=\frac{(-1)^{n-m+1}}{24}\Big(\sum_{l=1}^{m}\frac{1}{c_{l}}-\frac{n-2m+1}{2}\Big)\prod_{l=1}^{m}c_{l}^{c_{l}}.
\eea

We need to check that our conjectural formulae (\ref{58}) and (\ref{59}) match (\ref{62}) and (\ref{63}). First we give a lemma.
\begin{lemma}
If $n$ is odd, suppose $n=2r+1$, we have
\bea\label{64}
\Res_{w=0}\frac{(2w+1)\prod_{i=1}^{m}\prod_{s=1}^{c_{i}-1}(c_{i}w+s)}{(w+1)^{r-m+1}w^{r-m+1}}
=-\frac{1}{12}\Big(\frac{1}{c_{i}}-\sum_{i=1}^{m}c_{i}\Big)\cdot\prod_{i=1}^{m}c_{i}^{c_{i}-1}.
\eea
If $n$ is even, suppose $n=2r$, we have
\bea\label{65}
\Res_{w=0}\frac{\prod_{i=1}^{m}\prod_{s=1}^{c_{i}-1}(c_{i}w+s)}{(w+1)^{r-m}w^{r-m}}
=-\frac{1}{24}\Big(\frac{1}{c_{i}}-\sum_{i=1}^{m}c_{i}+\frac{3}{2}\Big)\cdot\prod_{i=1}^{m}c_{i}^{c_{i}-1}.
\eea
\end{lemma}
\emph{Proof}: The crucial point is to notice that
\ben
\Res_{w=0}\frac{(2w+1)\prod_{i=1}^{m}\prod_{s=1}^{c_{i}-1}(c_{i}w+s)}{(w+1)^{r-m+1}w^{r-m+1}}
=\Res_{w=-1}\frac{(2w+1)\prod_{i=1}^{m}\prod_{s=1}^{c_{i}-1}(c_{i}w+s)}{(w+1)^{r-m+1}w^{r-m+1}}
\een
and
\ben
\Res_{w=0}\frac{\prod_{i=1}^{m}\prod_{s=1}^{c_{i}-1}(c_{i}w+s)}{(w+1)^{r-m}w^{r-m}}
=\Res_{w=-1}\frac{\prod_{i=1}^{m}\prod_{s=1}^{c_{i}-1}(c_{i}w+s)}{(w+1)^{r-m}w^{r-m}}
\een
by  substitution of variables. Thus by the residue theorem on $\mathbb{P}^{1}$, it suffices to compute
\ben
\Res_{w=\infty}\frac{(2w+1)\prod_{i=1}^{m}\prod_{s=1}^{c_{i}-1}(c_{i}w+s)}{(w+1)^{r-m+1}w^{r-m+1}}
\een
and
\ben
\Res_{w=\infty}\frac{\prod_{i=1}^{m}\prod_{s=1}^{c_{i}-1}(c_{i}w+s)}{(w+1)^{r-m}w^{r-m}}.
\een
We leave the details to the reader.
\pqed

The functions $I_{0,q}(t)$ for $X$ are defined by
\bea\label{67}
\sum_{q=0}^{\infty}I_{0,q}(t)w^{q}=e^{wt}\sum_{d=0}^{\infty}e^{dt}\frac{\prod_{i=1}^{m}\prod_{s=0}^{c_{i}d-1}(-c_{i}w-s)}
{\prod_{s=1}^{d}(w+s)^{n-m+1}}=e^{wt}\Bigg[1+\sum_{d=1}^{\infty}e^{dt}\frac{(-1)^{(n-m+1)d}\prod_{i=1}^{m}\prod_{s=0}^{c_{i}d-1}(c_{i}w+s)}
{\prod_{s=1}^{d}(w+s)^{n-m+1}}\Bigg].\nn\\
\eea
For $q\geq p$, define
\bea\label{73}
I_{p,q}(t)=\frac{d}{dt}\Big(\frac{I_{p-1,q}(t)}{I_{p-1,p-1}(t)}\Big).
\eea
\begin{proposition}\label{66}
 For $2\nmid n$,
\bea\label{6}
\sum_{p=1}^{(n-1)/2}\frac{(n+1-2p)^{2}}{8}[e^{t}]\big(\ln I_{p,p}(t)\big)
=-\frac{(-1)^{\sum_{i=1}^{m}c_{i}}}{24}\Big(\sum_{i=1}^{m}\frac{1}{c_{i}}-\sum_{i=1}^{m}c_{i}\Big)\cdot\prod_{i=1}^{m}c_{i}^{c_{i}},
\eea
and for $2\mid n$,
\bea\label{8}
\sum_{p=1}^{(n-2)/2}\frac{(n+2-2p)(n-2p)}{8}[e^{t}]\big(\ln I_{p,p}(t)\big)
=-\frac{(-1)^{\sum_{i=1}^{m}c_{i}}}{24}\Big(\sum_{i=1}^{m}\frac{1}{c_{i}}-\sum_{i=1}^{m}c_{i}+\frac{3}{2}\Big)\cdot\prod_{i=1}^{m}c_{i}^{c_{i}}.
\eea
\end{proposition}

\emph{proof}: For a fixed $n$, suppose
\bea\label{11}
\frac{(-1)^{n-m+1}\prod_{i=1}^{m}\prod_{s=0}^{c_{i}-1}(c_{i}w+s)}{(w+1)^{n-m+1}}=a_{1}w+a_{2}w^{2}+\cdots,
\eea
then a straightforward induction shows
\bea\label{7}
I_{p,p}(t)=1+e^{t}\sum_{k=1}^{p}a_{k}\binom{p-1}{k-1}+O[e^{2t}].
\eea
Now we treat the cases that $n$ is odd or even separately.\\

(i)$n=2r+1$, and $r\geq 0$. By (\ref{7}) we have
\bea\label{14}
&&\sum_{p=1}^{(n-1)/2}\frac{(n+1-2p)^{2}}{8}[e^{t}]\big(\ln I_{p,p}(t)\big)=\sum_{p=1}^{r}\frac{(r+1-p)^{2}}{2}\sum_{k=1}^{p}a_{k}\binom{p-1}{k-1}\nn\\
&=&\sum_{k=1}^{r}a_{k}\sum_{p=1}^{r}\frac{(r+1-p)^{2}}{2}\binom{p-1}{k-1}.\nn\\
\eea
Since
\ben
(m+1-p)^{2}=(p+1)p-p(2m+3)+(m+1)^{2},
\een
we have
\bea\label{15}
&&\sum_{p=1}^{r}(r+1-p)^{2}\binom{p-1}{k-1}\nn\\&=&\sum_{p=1}^{r}
\Big((p+1)p\binom{p-1}{k-1}-(2r+3)p\binom{p-1}{k-1}+(r+1)^{2}\binom{p-1}{k-1}\Big)\nn\\
&=&\sum_{p=1}^{r}\Big((k+1)k\binom{p+1}{k+1}-(2r+3)k\binom{p}{k}+(r+1)^{2}\binom{p-1}{k-1}\Big)\nn\\
&=&(k+1)k\binom{r+2}{k+2}-(2r+3)k\binom{r+1}{k+1}+(r+1)^{2}\binom{r}{k}\nn\\
&=&\big((k+2)(k+1)-2(k+2)+2\big)\binom{r+2}{k+2}-(2r+3)\big((k+1)-1\big)\binom{r+1}{k+1}+(r+1)^{2}\binom{r}{k}\nn\\
&=&2\binom{r+2}{k+2}-\binom{r+1}{k+1}.\nn\\
\eea
Thus by (\ref{14}) and (\ref{15}) we have
\ben
&&\sum_{p=1}^{(n-1)/2}\frac{(n+1-2p)^{2}}{8}[e^{t}]\big(\ln I_{p,p}(t)\big)=
\frac{1}{2}\sum_{k=1}^{r}a_{k}\Bigg(2\binom{r+2}{k+2}-\binom{r+1}{k+1}\Bigg)\\
&=& \frac{1}{2}[w^{-2}]\Bigg(\frac{2(-1)^{n-m+1}\prod_{i=1}^{m}\prod_{s=0}^{c_{i}-1}(c_{i}w+s)}{(w+1)^{n-m+1}}
\cdot \Big(1+\frac{1}{w}\Big)^{r+2}\Bigg)\\
&&-\frac{1}{2}[w^{-1}]\Bigg(\frac{(-1)^{n-m+1}\prod_{i=1}^{m}\prod_{s=0}^{c_{i}-1}(c_{i}w+s)}{(w+1)^{n-m+1}}
\cdot \Big(1+\frac{1}{w}\Big)^{r+1}\Bigg)\\
&=&\frac{(-1)^{n-m+1}\prod_{i=1}^{m}c_{i}}{2}[w^{r-m}]\Bigg(\frac{(2w+1)\prod_{i=1}^{m}\prod_{s=1}^{c_{i}-1}(c_{i}w+s)}{(w+1)^{r-m+1}}
\Bigg).
\een
Thus by (\ref{64}) we obtain (\ref{6}).\\

(ii)$n=2r$, and $r\geq 1$. By (\ref{7}) we have
\bea\label{9}
&&\sum_{p=1}^{(n-2)/2}\frac{(n+2-2p)(n-2p)}{8}[e^{t}]\big(\ln I_{p,p}(t)\big)=\sum_{p=1}^{r-1}\frac{(r+1-p)(r-p)}{2}\sum_{k=1}^{p}a_{k}\binom{p-1}{k-1}\nn\\
&=&\sum_{k=1}^{r-1}a_{k}\sum_{p=1}^{r-1}\frac{(r+1-p)(r-p)}{2}\binom{p-1}{k-1}.
\eea
A similar computation as in the $n$ odd case shows
\bea\label{10}
\sum_{p=1}^{r-1}(r+1-p)(r-p)\binom{p-1}{k-1}
=2\binom{r+1}{k+2}.
\eea
Thus by (\ref{9}) and (\ref{10}), we see
\ben
&&\sum_{p=1}^{(n-2)/2}\frac{(n-2p)(n-2-2p)}{8}[e^{t}]\big(\ln I_{p,p}(t)\big)=
\sum_{k=1}^{r}a_{k}\binom{r+1}{k+2}\\
&=&[w^{-2}]\Bigg(\frac{(-1)^{n-m+1}\prod_{i=1}^{m}\prod_{s=0}^{c_{i}-1}(c_{i}w+s)}{(w+1)^{n-m+1}}
\cdot \Big(1+\frac{1}{w}\Big)^{r+1}\Bigg)\\
&=&(-1)^{n-m+1}\prod_{i=1}^{m}c_{i}\cdot[w^{r-m-1}]\Bigg(\frac{(2w+1)\prod_{i=1}^{m}\prod_{s=1}^{c_{i}-1}(c_{i}w+s)}{(w+1)^{r-m}}
\Bigg).
\een
Then (\ref{8}) follows from (\ref{65}).
\pqed

When $m\geq 2$, Prop.\ref{66} together with the contribution from $-\frac{n+1}{48}\ln(1-\prod_{i=1}^{m}(-c_{i})^{c_{i}}e^{t})$ or $-\frac{n-2}{48}\ln(1-\prod_{i=1}^{m}(-c_{i})^{c_{i}}e^{t})$ ($n$ is odd or even, resp.) gives (\ref{63}). When $m=1$, from (\ref{67}) it is easy to see
\ben
T=t+\sum_{d=1}^{\infty}e^{dt}\frac{(-1)^{nd}}{d}\frac{(nd)!}{(d!)^{n}}.
\een
Take this into account, we also recover  (\ref{62}). So we have proved
\begin{theorem}
The conjecture \ref{60} holds for all degree one invariants.
\end{theorem}
\pqed

\begin{remark}
The same method shows that for the Calabi-Yau hypersurface $Y$ in $\mathbb{P}^{n-1}$ we have
\ben
N_{1,1}^{Y}=n!\Bigg[\big(\frac{(n-2)(n+1)}{48}+\frac{1-(1-n)^{n}}{24n^{2}}\big)\sum_{s=2}^{n}\frac{n}{s}
+\frac{n^{2}-1+(1-n)^{n}}{24n}\Bigg]
-\frac{n^{n-1}(n-1)(n+2)}{48}.
\een
\end{remark}
\section{Two extremal cases}
In general as the degree $d$ increase, the graphs and their contributions  corresponding to the fixed loci will become more and more complicated, and thus a direct computation through virtual torus localization seems very difficult. But  for some special target spaces we can make a good choice of the linearization so that a lot of graphs give zero contributions (see, e.g., \cite{GP}). In principle, the larger $m$ is, the more flexible the choice of the linearizations is. We shall consider the two extremal cases: $X=\mathrm{Tot}\big(\mathcal{O}(-1)^{\oplus (l+1)}\rightarrow \mathbb{P}^{l}\big)$ and $X=\mathrm{Tot}\big(\mathcal{O}(-1)^{\oplus(l-1)}\oplus\mathcal{O}(-2)\rightarrow \mathbb{P}^{l}\big)$. In these two cases it is easy to see from (\ref{67}) and (\ref{73}) that $I_{p,p}(t)=1$ for $p$ in the ranges that appear in (\ref{58}) and (\ref{59}). So to prove conjecture \ref{60} in these two cases is equivalent to show
\begin{theorem}
For $X=\mathrm{Tot}\big(\mathcal{O}(-1)^{\oplus (l+1)}\rightarrow \mathbb{P}^{l}\big)$ we have
\bea
N_{1,d}=\frac{(-1)^{(l-1)d}(l+1)}{24d}.
\eea
For $X=\mathrm{Tot}\big(\mathcal{O}(-1)^{\oplus(l-1)}\oplus\mathcal{O}(-2)\rightarrow \mathbb{P}^{l}\big)$ we have
\bea
N_{1,d}^{X}=
\frac{(-1)^{(l-1)d}(l-1)4^{d}}{24d}.
\eea
\end{theorem}
In the following we treat the two cases separately. The choice of linearizations are following those of the similar cases in \cite{KP}
and \cite{PZ}. In the following computations we shall make repeatedly use of $\lambda_{1}^{2}=0$ on $\Mbar_{1,m}$ for $m\geq 1$, for example from this we have $\Lambda_{1}^{\vee}(x)\Lambda_{1}^{\vee}(-x)=-x^{2}$.
\subsubsection{$\mathcal{O}(-1)^{\oplus (l+1)}\rightarrow \mathbb{P}^{l}$}
Write $\mathcal{O}(-1)^{\oplus (l+1)}=\bigoplus_{i=1}^{l+1}L_{i}$, and choose torus linearizations on $L_{i}$ with weight $\alpha_{i}-\alpha_{k}$ at $P_{k}$, for $1\leq i,k\leq l+1$. In particular, $L_{i}$ has weight zero at $P_{i}$. The fixed loci with nonzero contributions are of the form
$$ \Gamma_{ij}=\xy
(0,0); (10,0), **@{-};
(0,0)*+{\circ};(10,0)*+{\bullet};(5,2)*+{d};
(0,3)*+{i};(10,3)*+{j};
\endxy,
$$
where $1\leq i\neq j\leq l+1$.
The contribution of $\Gamma_{ij}$ is
\ben
&&\frac{1}{d}\int_{\Mbar_{1,1}}\frac{\frac{\alpha_{j}-\alpha_{i}}{d}\prod_{k\neq i}\Lambda_{1}^{\vee}(\alpha_{i}-\alpha_{k})\cdot \prod_{k=1}^{l+1}\Big(\Lambda_{1}^{\vee}(\alpha_{k}-\alpha_{i}) \prod_{a=1}^{d-1}(\alpha_{k}-\alpha_{j}+a\frac{\alpha_{j}-\alpha_{i}}{d})\Big)}
{(\frac{\alpha_{i}-\alpha_{j}}{d}-\psi)
(\frac{d!}{d^{d}})^{2}(\alpha_{i}-\alpha_{j})^d(\alpha_{j}-\alpha_{i})^d\prod_{k\neq i,j}\prod_{a=0}^{d}(\alpha_{i}-\alpha_{k}+a\frac{\alpha_{j}-\alpha_{i}}{d})}\\
&=&-\frac{1}{24d}\frac{(\alpha_{j}-\alpha_{i})\prod_{k\neq i}(\alpha_{i}-\alpha_{k})(\alpha_{k}-\alpha_{i})\cdot \prod_{k=1}^{l+1} \prod_{a=1}^{d-1}(\alpha_{k}-\alpha_{j}+a\frac{\alpha_{j}-\alpha_{i}}{d})}
{(\alpha_{i}-\alpha_{j})
(\frac{d!}{d^{d}})^{2}(\alpha_{i}-\alpha_{j})^d(\alpha_{j}-\alpha_{i})^d\prod_{k\neq i,j}\prod_{a=0}^{d}(\alpha_{i}-\alpha_{k}+a\frac{\alpha_{j}-\alpha_{i}}{d})}.
\een
Note that
\ben
\alpha_{i}-\alpha_{k}+(d-a)\frac{\alpha_{j}-\alpha_{i}}{d}=-\Big(\alpha_{k}-\alpha_{j}+a\frac{\alpha_{j}-\alpha_{i}}{d}\Big),
\een
so the contribution is
\ben
\frac{(-1)^{(l-1)d}}{24d}\frac{\prod_{k\neq i,j}(\alpha_{i}-\alpha_{k})}{\prod_{k\neq i,j}(\alpha_{j}-\alpha_{k})}.
\een
Since  for any fixed $i$ we have
\ben
\sum_{j\neq i}\frac{\prod_{k\neq i,j}(\alpha_{i}-\alpha_{k})}{\prod_{k\neq i,j}(\alpha_{j}-\alpha_{k})}=1,
\een
we obtain
\ben
N_{1,d}=\frac{(-1)^{(l-1)d}(l+1)}{24d}.
\een

\subsubsection{ $\mathcal{O}(-1)^{\oplus(l-1)}\oplus\mathcal{O}(-2)\rightarrow \mathbb{P}^{l}$}
Choose the linearizations on $L_{i}$ such that for $1\leq i \leq l-1$, $L_{i}$ has weight $\alpha_{i}-\alpha_{k}$ at $P_{k}$, and $L_{l}$ has weight $\alpha_{l}+\alpha_{l+1}-2\alpha_{k}$ at $P_{k}$, $1\leq k\leq l+1$.
 The fixed loci which may have nonzero contributions are of three types.\\

 Type I:

$$ \Gamma_{s;k_{1},d_{1};\cdots;k_{m},d_{m}}=\xy
(0,0); (10,10), **@{-};(0,0);(10,8), **@{-};(0,0);(10,-8), **@{-};(0,0);(10,-10), **@{-};
(0,0)*+{\circ};(10,10)*+{\bullet};(10,8)*+{\bullet};(10,-8)*+{\bullet};(10,-10)*+{\bullet};(10,6)*+{\cdot};(10,4)*+{\cdot};(10,2)*+{\cdot};
(10,-2)*+{\cdot};(10,-4)*+{\cdot};(10,-6)*+{\cdot};
(-1,3)*+{s};(13,10)*+{k_{1}}; (13,7.5)*+{k_{2}};(13,-10)*+{k_{m}};
\endxy,
$$
where $1\leq k_{1},\cdots,k_{m}\leq l-1$, $m\geq 1$, with edges of degree $d_{1},\cdots,d_{m}$ respectively, and $s=l$ or $l+1$. When $s=l+1$, the contribution is
\ben
&&\frac{1}{|\mathrm{Aut}(\Gamma_{s;k_{1},\cdots,k_{m}})|\prod_{i=1}^{m}d_{i}}
\int_{\Mbar_{1,m}}\prod_{i=1}^{m}\frac{\alpha_{k_{i}}-\alpha_{l+1}}{d_{i}}
\prod_{j=1}^{l}\Big((\alpha_{l+1}-\alpha_{j})^{m-1}\Lambda_{1}^{\vee}(\alpha_{l+1}-\alpha_{j})\Big)\\
&&\prod_{j=1}^{l-1}\Big((\alpha_{j}-\alpha_{l+1})^{m-1}\Lambda_{1}^{\vee}(\alpha_{j}-\alpha_{l+1})
\prod_{i=1}^{m}\prod_{a=1}^{d_{i}-1}(\alpha_{j}-\alpha_{l+1}+a\frac{\alpha_{l+1}-\alpha_{k_{i}}}{d_{i}})\Big)\\
&&\frac{(\alpha_{l}-\alpha_{l+1})^{m-1}\Lambda_{1}^{\vee}(\alpha_{l}-\alpha_{l+1})\prod_{i=1}^{m}
\prod_{a=1}^{2d_{i}-1}(\alpha_{l}-\alpha_{l+1}+a\frac{\alpha_{l+1}-\alpha_{k_{i}}}{d_{i}})}
{\prod_{i=1}^{m}\Big((\frac{\alpha_{l+1}-\alpha_{k_{i}}}{d_{i}}-\psi_{i})\frac{(d_{i}!)^{2}}{d_{i}^{2d_{i}}}(-1)^{d_{i}}
(\alpha_{l+1}-\alpha_{k_{i}})^{2d_{i}}
\prod_{r\neq l+1,k_{i}}\prod_{a=0}^{d_{i}}(\alpha_{l+1}-\alpha_{r}+a\frac{\alpha_{k_{i}}-\alpha_{l+1}}{d_{i}})\Big)}\\
&=&\frac{1}{|\mathrm{Aut}(\Gamma_{k_{1},\cdots,k_{m}})|\prod_{i=1}^{m}d_{i}^{2}}
\int_{\Mbar_{1,m}}
\frac{1}{\prod_{i=1}^{m}(\frac{\alpha_{l+1}-\alpha_{k_{i}}}{d_{i}}-\psi_{i})}\\
&&\frac{(-1)^{(l-1)d+m}(\alpha_{l+1}-\alpha_{l})^{m}\prod_{j=1}^{l-1}(\alpha_{l+1}-\alpha_{j})^{m}
\prod_{i=1}^{m}\Big(
\prod_{a=0}^{d_{i}-1}(\alpha_{l}-\alpha_{k_{i}}+a\frac{\alpha_{l+1}-\alpha_{k_{i}}}{d_{i}})\Big)}
{\prod_{i=1}^{m}\Big(\frac{d_{i}!}{d_{i}^{d_{i}}}
(\alpha_{l+1}-\alpha_{k_{i}})^{d_{i}-1}(\alpha_{k_{i}}-\alpha_{l})
\prod_{r\neq l+1,l,k_{i}}(\alpha_{k_{i}}-\alpha_{r})
\Big)}.
\een

Similarly  When $s=l+1$, the contribution is
\ben
&&\frac{1}{|\mathrm{Aut}(\Gamma_{k_{1},\cdots,k_{m}})|\prod_{i=1}^{m}d_{i}^{2}}
\int_{\Mbar_{1,m}}
\frac{1}{\prod_{i=1}^{m}(\frac{\alpha_{l+1}-\alpha_{k_{i}}}{d_{i}}-\psi_{i})}\\
&&\frac{(-1)^{(l-1)d+m}(\alpha_{l+1}-\alpha_{l})^{m}\prod_{j=1}^{l-1}(\alpha_{l+1}-\alpha_{j})^{m}
\prod_{i=1}^{m}\Big(
\prod_{a=0}^{d_{i}-1}(\alpha_{l}-\alpha_{k_{i}}+a\frac{\alpha_{l+1}-\alpha_{k_{i}}}{d_{i}})\Big)}
{\prod_{i=1}^{m}\Big(\frac{d_{i}!}{d_{i}^{d_{i}}}
(\alpha_{l+1}-\alpha_{k_{i}})^{d_{i}-1}(\alpha_{k_{i}}-\alpha_{l})
\prod_{r\neq l+1,l,k_{i}}(\alpha_{k_{i}}-\alpha_{r})
\Big)}.
\een
The crucial observation is that, in these contributions the factor $\alpha_{l+1}-\alpha_{l}$ appears at least once, and we shall see that $\alpha_{l+1}-\alpha_{l}$ does not appear in the denominator of the sums of the contributions of the other types. So we are able to set $\alpha_{l+1}=\alpha_{l}$ and thus the type I graphs contribute nothing.\\

Type II:

$$ \Gamma_{s;k_{0},d_{0};k_{1},d_{1};\cdots;k_{m},d_{m}}^{\prime}=
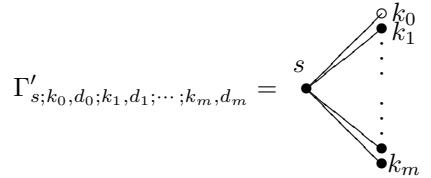
\xy
(0,0); (10,10), **@{-};(0,0);(10,8), **@{-};(0,0);(10,-8), **@{-};(0,0);(10,-10), **@{-};
(0,0)*+{\bullet};(10,10)*+{\circ};(10,8)*+{\bullet};(10,-8)*+{\bullet};(10,-10)*+{\bullet};(10,6)*+{\cdot};(10,4)*+{\cdot};(10,2)*+{\cdot};
(10,-2)*+{\cdot};(10,-4)*+{\cdot};(10,-6)*+{\cdot};
(-1,3)*+{s};(13,10)*+{k_{0}}; (13,7.5)*+{k_{1}};(13,-10)*+{k_{m}};
\endxy,
$$
where $1\leq k_{0}, k_{1},\cdots,k_{m}\leq l-1$, $m\geq 0$, with edges of degree $d_{0}, d_{1},\cdots,d_{m}$ respectively, and $s=l$ or $l+1$.
When $s=l+1$, the contribution is
\ben
&&\frac{1}{|\mathrm{Aut}(\Gamma_{s;k_{0},k_{1},\cdots,k_{m}}^{\prime})|\prod_{i=0}^{m}d_{i}}
\int_{\Mbar_{1,1}}\frac{\prod_{j\neq k_{0}}\Lambda_{1}^{\vee}(\alpha_{k_{0}}-\alpha_{j})
\prod_{i=1}^{l-1}\Lambda_{1}^{\vee}(\alpha_{i}-\alpha_{k_{0}})\cdot\Lambda_{1}^{\vee}(\alpha_{l}+\alpha_{l+1}-2\alpha_{k_{0}})}
{\frac{\alpha_{k_{0}}-\alpha_{l+1}}{d_{0}}-\psi}\\
&&\int_{\Mbar_{0,m+1}}\frac{1}{\prod_{i=0}^{m}(\frac{\alpha_{l+1}-\alpha_{k_{i}}}{d_{i}}-\psi_{i})}\cdot
\prod_{i=1}^{m}\frac{\alpha_{k_{i}}-\alpha_{l+1}}{d_{i}}
\prod_{j=1}^{l}\Big((\alpha_{l+1}-\alpha_{j})^{m}\Big)\\
&&\prod_{j=1}^{l-1}\Big((\alpha_{j}-\alpha_{l+1})^{m}
\prod_{i=0}^{m}\prod_{a=1}^{d_{i}-1}(\alpha_{j}-\alpha_{l+1}+a\frac{\alpha_{l+1}-\alpha_{k_{i}}}{d_{i}})\Big)\\
&&\frac{(\alpha_{l}-\alpha_{l+1})^{m}\prod_{i=0}^{m}
\prod_{a=1}^{2d_{i}-1}(\alpha_{l}-\alpha_{l+1}+a\frac{\alpha_{l+1}-\alpha_{k_{i}}}{d_{i}})}
{\prod_{i=0}^{m}\Big(\frac{(d_{i}!)^{2}}{d_{i}^{2d_{i}}}(-1)^{d_{i}}
(\alpha_{l+1}-\alpha_{k_{i}})^{2d_{i}}
\prod_{r\neq l+1,k_{i}}\prod_{a=0}^{d_{i}}(\alpha_{l+1}-\alpha_{r}+a\frac{\alpha_{k_{i}}-\alpha_{l+1}}{d_{i}})\Big)}\\
&=&\frac{d_{0}}{|\mathrm{Aut}(\Gamma_{s;k_{0},k_{1},\cdots,k_{m}}^{\prime})|\prod_{i=0}^{m}d_{i}^{2}}
\Big(-\frac{1}{24}(-1)^{l-2}d_{0}(\alpha_{k_{0}}-\alpha_{l})\prod_{j=1,\neq k_{0}}^{l-1}(\alpha_{k_{0}}-\alpha_{j})^{2}
(\alpha_{l}+\alpha_{l+1}-2\alpha_{k_{0}})\Big)\\
&&\int_{\Mbar_{0,m+1}}\frac{1}{\prod_{i=0}^{m}(\frac{\alpha_{l+1}-\alpha_{k_{i}}}{d_{i}}-\psi_{i})}\cdot
\frac{1}{\alpha_{k_{0}}-\alpha_{l+1}}\cdot(-1)^{lm}
\prod_{j=1}^{l}\Big((\alpha_{l+1}-\alpha_{j})^{2m}\Big)\\
&&\frac{\prod_{i=0}^{m}\Big((-1)^{(l-1)(d_{i}-1)}
\prod_{a=d_{i}}^{2d_{i}-1}(\alpha_{l}-\alpha_{l+1}+a\frac{\alpha_{l+1}-\alpha_{k_{i}}}{d_{i}})\Big)}
{\prod_{i=0}^{m}\Big(\frac{d_{i}!}{d_{i}^{d_{i}}}
(\alpha_{l+1}-\alpha_{k_{i}})^{d_{i}}
\prod_{r\neq l+1,l,k_{i}}(\alpha_{l+1}-\alpha_{r})(\alpha_{k_{i}}-\alpha_{r})\cdot(\alpha_{l+1}-\alpha_{l})(\alpha_{k_{i}}-\alpha_{l})
\Big)}.
\een
When $m>0$, the power of $\alpha_{l+1}-\alpha_{l}$ in the numerator is not less than that in the denominator. To show that the sums of contributions of the type II graphs has no factor of $\alpha_{l+1}-\alpha_{l}$ in its denominator, we only need to consider the $m=0$ case. When $m=0$, the above contribution is
\ben
\frac{(-1)^{(l-1)d_{0}+1}d_{0}^{d_{0}-1}}{24d_{0}!}\frac{\prod_{j=1,\neq k_{0}}^{l-1}(\alpha_{k_{0}}-\alpha_{j})\cdot
\prod_{a=0}^{d_{0}}(\alpha_{l}-\alpha_{k_{0}}+a\frac{\alpha_{l+1}-\alpha_{k_{0}}}{d_{0}})}
{(\alpha_{l+1}-\alpha_{k_{0}})^{d_{0}}
\prod_{r\neq l+1,l,k_{0}}(\alpha_{l+1}-\alpha_{r})\cdot(\alpha_{l+1}-\alpha_{l})}.
\een
Thus the sum of the contributions of $\Gamma_{l;k_{0}}^{\prime}$ and $\Gamma_{l+1;k_{0}}^{\prime}$ is
\ben
&&\frac{(-1)^{(l-1)d+1}d^{d-1}}{24d!}\prod_{j=1,\neq k_{0}}^{l-1}(\alpha_{k_{0}}-\alpha_{j})\\
&&\cdot\frac{1}{\alpha_{l+1}-\alpha_{l}}\Big[\frac{
\prod_{a=0}^{d}(\alpha_{l}-\alpha_{k_{0}}+a\frac{\alpha_{l+1}-\alpha_{k_{0}}}{d})}
{(\alpha_{l+1}-\alpha_{k_{0}})^{d}
\prod_{r\neq l+1,l,k_{0}}(\alpha_{l+1}-\alpha_{r})}
-\frac{
\prod_{a=0}^{d}(\alpha_{l+1}-\alpha_{k_{0}}+a\frac{\alpha_{l}-\alpha_{k_{0}}}{d})}
{(\alpha_{l}-\alpha_{k_{0}})^{d}
\prod_{r\neq l+1,l,k_{0}}(\alpha_{l}-\alpha_{r})}\Big]\\
&=&\frac{(-1)^{(l-1)d+1}d^{d-1}}{24d!}\frac{\prod_{j=1,\neq k_{0}}^{l-1}(\alpha_{k_{0}}-\alpha_{j})}
{(\alpha_{l+1}-\alpha_{k_{0}})^{d}(\alpha_{l}-\alpha_{k_{0}})^{d}
\cdot\prod_{r\neq l+1,l,k_{0}}(\alpha_{l+1}-\alpha_{r})(\alpha_{l}-\alpha_{r})}\\
&&\cdot\frac{1}{{\alpha_{l+1}-\alpha_{l}}}
\cdot\Big[\prod_{a=0}^{d}(\alpha_{l}-\alpha_{k_{0}}+a\frac{\alpha_{l+1}-\alpha_{k_{0}}}{d})
\cdot(\alpha_{l}-
\alpha_{k_{0}})^{d}\prod_{r\neq l+1,l,k_{0}}(\alpha_{l}-\alpha_{r})\\
&&
-\prod_{a=0}^{d}(\alpha_{l+1}-\alpha_{k_{0}}+a\frac{\alpha_{l}-\alpha_{k_{0}}}{d})\cdot
(\alpha_{l+1}-\alpha_{k_{0}})^{d}
\prod_{r\neq l+1,l,k_{0}}(\alpha_{l+1}-\alpha_{r})\Big].
\een
The sum of the group of terms in the square brackets of the last expression is divisible by $\alpha_{l+1}-\alpha_{l}$. Therefore we have shown that the sum of the contributions of type II graphs has no factor $\alpha_{l+1}-\alpha_{l}$ in its denominator. We shall see the type III contribution also has no factor $\alpha_{l+1}-\alpha_{l}$ in the denominator. So we are able to set $\alpha_{l+1}=\alpha_{l}$. Then we see that
a type II graph has no contribution unless $m=0$ or $m=1$. Now we compute the contributions of $m=0$ and $m=1$ cases separately. Since
\ben
&&\prod_{a=0}^{d}(\alpha_{l}-\alpha_{k_{0}}+a\frac{\alpha_{l+1}-\alpha_{k_{0}}}{d})\\
&=&(\alpha_{l}-\alpha_{k_{0}})\frac{d!}{d^{d}}\prod_{a=1}^{d}(\frac{d(\alpha_{l}-\alpha_{k_{0}})}{a}+\alpha_{l+1}-\alpha_{k_{0}})\\
&=&(\alpha_{l}-\alpha_{k_{0}})\frac{d!}{d^{d}}\Big(
\frac{(2d)!}{(d!)^{2}}(\alpha_{l}-\alpha_{k_{0}})^{d}
+(\alpha_{l+1}-\alpha_{l})\frac{(2d)!}{(d!)^{2}}(\alpha_{l}-\alpha_{k_{0}})^{d-1}\sum_{a=1}^{d}\frac{a}{d+a}+O[(\alpha_{l+1}-\alpha_{l})^{2}]\Big)\\
&=&
\frac{(2d)!}{d!d^{d}}(\alpha_{l}-\alpha_{k_{0}})^{d+1}
+(\alpha_{l+1}-\alpha_{l})\frac{(2d)!}{d!d^{d}}(\alpha_{l}-\alpha_{k_{0}})^{d}\sum_{a=1}^{d}\frac{a}{d+a}+O[(\alpha_{l+1}-\alpha_{l})^{2}],
\een

\ben
&&\prod_{a=0}^{d}(\alpha_{l+1}-\alpha_{k_{0}}+a\frac{\alpha_{l}-\alpha_{k_{0}}}{d})\\
&=&\frac{(2d)!}{d!d^{d}}(\alpha_{l}-\alpha_{k_{0}})^{d+1}+
(\alpha_{l+1}-\alpha_{l})\frac{(2d)!}{d!d^{d}}(\alpha_{l}-\alpha_{k_{0}})^{d}\sum_{a=0}^{d}\frac{d}{d+a}
+O[(\alpha_{l+1}-\alpha_{l})^{2}],
\een

\ben
(\alpha_{l+1}-\alpha_{k_{0}})^{d}
=(\alpha_{l}-\alpha_{k_{0}})^{d}+d(\alpha_{l+1}-\alpha_{l})(\alpha_{l}-\alpha_{k_{0}})^{d-1}+O[(\alpha_{l+1}-\alpha_{l})^{2}],
\een

\ben
&&\prod_{r\neq l+1,l,k_{0}}(\alpha_{l+1}-\alpha_{r})\\
&=&\prod_{r\neq l+1,l,k_{0}}(\alpha_{l}-\alpha_{r})
+(\alpha_{l+1}-\alpha_{l})\prod_{r\neq l+1,l,k_{0}}(\alpha_{l}-\alpha_{r})\sum_{r\neq l+1,l,k_{0}}\frac{1}{\alpha_{l}-\alpha_{r}}
+O[(\alpha_{l+1}-\alpha_{l})^{2}],
\een

we have
\ben
&&\prod_{a=0}^{d}(\alpha_{l}-\alpha_{k_{0}}+a\frac{\alpha_{l+1}-\alpha_{k_{0}}}{d})
\cdot(\alpha_{l}-
\alpha_{k_{0}})^{d}\prod_{r\neq l+1,l,k_{0}}(\alpha_{l}-\alpha_{r})\\
&&
-\prod_{a=0}^{d}(\alpha_{l+1}-\alpha_{k_{0}}+a\frac{\alpha_{l}-\alpha_{k_{0}}}{d})\cdot
(\alpha_{l+1}-\alpha_{k_{0}})^{d}
\prod_{r\neq l+1,l,k_{0}}(\alpha_{l+1}-\alpha_{r})\\
&=&(\alpha_{l+1}-\alpha_{l})\frac{(2d)!}{d!d^{d}}(\alpha_{l}-\alpha_{k_{0}})^{d}\sum_{a=1}^{d}\frac{a}{d+a}
\cdot(\alpha_{l}-
\alpha_{k_{0}})^{d}\prod_{r\neq l+1,l,k_{0}}(\alpha_{l}-\alpha_{r})\\
&&-(\alpha_{l+1}-\alpha_{l})\frac{(2d)!}{d!d^{d}}(\alpha_{l}-\alpha_{k_{0}})^{d}\sum_{a=0}^{d}\frac{d}{d+a}
\cdot(\alpha_{l}-
\alpha_{k_{0}})^{d}\prod_{r\neq l+1,l,k_{0}}(\alpha_{l}-\alpha_{r})\\
&&-d(\alpha_{l+1}-\alpha_{l})(\alpha_{l}-\alpha_{k_{0}})^{d-1}\cdot\frac{(2d)!}{d!d^{d}}(\alpha_{l}-\alpha_{k_{0}})^{d+1}
\prod_{r\neq l+1,l,k_{0}}(\alpha_{l}-\alpha_{r})\\
&&-(\alpha_{l+1}-\alpha_{l})\prod_{r\neq l+1,l,k_{0}}(\alpha_{l}-\alpha_{r})\sum_{r\neq l+1,l,k_{0}}\frac{1}{\alpha_{l}-\alpha_{r}}\cdot
\frac{(2d)!}{d!d^{d}}(\alpha_{l}-\alpha_{k_{0}})^{d+1}(\alpha_{l}-\alpha_{k_{0}})^{d}+O[(\alpha_{l+1}-\alpha_{l})^{2}]\\
&=&(\alpha_{l+1}-\alpha_{l})\frac{(2d)!}{d!d^{d}}(\alpha_{l}-\alpha_{k_{0}})^{2d}\prod_{r\neq l+1,l,k_{0}}(\alpha_{l}-\alpha_{r})\\
&&\cdot\Big(2\sum_{a=1}^{d}\frac{a}{d+a}-2d-1-
(\alpha_{l}-\alpha_{k_{0}})\sum_{r\neq l+1,l,k_{0}}\frac{1}{\alpha_{l}-\alpha_{r}}\Big)+O[(\alpha_{l+1}-\alpha_{l})^{2}].
\een
Thus setting $\alpha_{l+1}=\alpha_{l}=\alpha$, the sum of the contributions of $\Gamma_{l;k_{0}}^{\prime}$ and $\Gamma_{l+1;k_{0}}^{\prime}$ is
\ben
&&\frac{(-1)^{(l-1)d+1}d^{d-1}}{24d!}\frac{\prod_{j=1,\neq k_{0}}^{l-1}(\alpha_{k_{0}}-\alpha_{j})}
{(\alpha_{l+1}-\alpha_{k_{0}})^{d}(\alpha_{l}-\alpha_{k_{0}})^{d}
\cdot\prod_{r\neq l+1,l,k_{0}}(\alpha_{l+1}-\alpha_{r})(\alpha_{l}-\alpha_{r})}\cdot\\
&&\frac{(2d)!}{d!d^{d}}(\alpha_{l}-\alpha_{k_{0}})^{2d}\prod_{r\neq l+1,l,k_{0}}(\alpha_{l}-\alpha_{r})
\cdot\Big(2\sum_{a=1}^{d}\frac{a}{d+a}-2d-1-
(\alpha_{l}-\alpha_{k_{0}})\sum_{r\neq l+1,l,k_{0}}\frac{1}{\alpha_{l}-\alpha_{r}}\Big)\\
&=&\frac{(-1)^{(l-1)d}(2d)!}{24(d!)^{2}d}\frac{\prod_{j=1,\neq k_{0}}^{l-1}(\alpha_{k_{0}}-\alpha_{j})}
{\cdot\prod_{j=1,\neq k_{0}}^{l-1}(\alpha-\alpha_{j})}
\Big(2d\sum_{a=1}^{d}\frac{1}{d+a}+1+
(\alpha-\alpha_{k_{0}})\sum_{j=1,\neq k_{0}}^{l-1}\frac{1}{\alpha-\alpha_{j}}\Big).
\een

The contribution of $\Gamma_{l+1;k_{0},k_{1}}^{\prime}$ is
\ben
&&\frac{1}{d_{0}d_{1}^{2}}
\Big(-\frac{1}{24}(-1)^{l-2}d_{0}(\alpha_{k_{0}}-\alpha_{l})\prod_{j=1,\neq k_{0}}^{l-1}(\alpha_{k_{0}}-\alpha_{j})^{2}\cdot
(\alpha_{l}+\alpha_{l+1}-2\alpha_{k_{0}})\Big)\\
&&\int_{\Mbar_{0,2}}\frac{1}{\prod_{i=0}^{1}(\frac{\alpha_{l+1}-\alpha_{k_{i}}}{d_{i}}-\psi_{i})}\cdot
\frac{1}{\alpha_{k_{0}}-\alpha_{l+1}}\cdot(-1)^{l}
\prod_{j=1}^{l}\Big((\alpha_{l+1}-\alpha_{j})^{2}\Big)\\
&&\frac{\prod_{i=0}^{1}\Big((-1)^{(l-1)(d_{i}-1)}
\prod_{a=d_{i}}^{2d_{i}-1}(\alpha_{l}-\alpha_{l+1}+a\frac{\alpha_{l+1}-\alpha_{k_{i}}}{d_{i}})\Big)}
{\prod_{i=0}^{1}\Big(\frac{d_{i}!}{d_{i}^{d_{i}}}
(\alpha_{l+1}-\alpha_{k_{i}})^{d_{i}}
\prod_{r\neq l+1,l,k_{i}}(\alpha_{l+1}-\alpha_{r})(\alpha_{k_{i}}-\alpha_{r})\cdot(\alpha_{l+1}-\alpha_{l})(\alpha_{k_{i}}-\alpha_{l})
\Big)}\\
&=&\frac{(-1)^{(l-1)d}d_{0}^{d_{0}+1}d_{1}^{d_{1}-1}}{24d_{0}!d_{1}!}
\Big(\prod_{j=1,\neq k_{0}}^{l-1}(\alpha_{k_{0}}-\alpha_{j})\cdot
(\alpha_{l}+\alpha_{l+1}-2\alpha_{k_{0}})\Big)\\
&&\frac{1}{d\alpha_{l+1}-d_{1}\alpha_{k_{0}}-d_{0}\alpha_{k_{1}}}\cdot
(\alpha_{l+1}-\alpha_{k_{1}})\\
&&\frac{\prod_{i=0}^{1}\Big(
\prod_{a=0}^{d_{i}-1}(\alpha_{l}-\alpha_{k_{i}}+a\frac{\alpha_{l+1}-\alpha_{k_{i}}}{d_{i}})\Big)}
{(\alpha_{l+1}-\alpha_{k_{0}})^{d_{0}}
(\alpha_{l+1}-\alpha_{k_{1}})^{d_{1}}
\prod_{j=1,\neq k_{1}}^{l-1}(\alpha_{k_{1}}-\alpha_{j})\cdot(\alpha_{k_{1}}-\alpha_{l})},
\een
setting $\alpha_{l+1}=\alpha_{l}=\alpha$, the above contribution becomes
\ben
&&\frac{(-1)^{(l-1)d+1}d_{0}(2d_{0})!(2d_{1})!}{48(d_{0}!)^{2}(d_{1}!)^{2}d_{1}}
\prod_{j=1,\neq k_{0}}^{l-1}(\alpha_{k_{0}}-\alpha_{j})
\cdot
\frac{1}
{\prod_{j=1,\neq k_{1}}^{l-1}(\alpha_{k_{1}}-\alpha_{j})}\cdot
\frac{\alpha-\alpha_{k_{0}}}{d\alpha-d_{1}\alpha_{k_{0}}-d_{0}\alpha_{k_{1}}}.
\een
Therefore the sum of the contributions of $\Gamma_{l+1;k_{0},k_{1}}^{\prime}$ and $\Gamma_{l;k_{0},k_{1}}^{\prime}$ is
\ben
\frac{(-1)^{(l-1)d+1}d_{0}(2d_{0})!(2d_{1})!}{24(d_{0}!)^{2}(d_{1}!)^{2}d_{1}}
\prod_{j=1,\neq k_{0}}^{l-1}(\alpha_{k_{0}}-\alpha_{j})
\cdot
\frac{1}
{\prod_{j=1,\neq k_{1}}^{l-1}(\alpha_{k_{1}}-\alpha_{j})}\cdot
\frac{\alpha-\alpha_{k_{0}}}{d\alpha-d_{1}\alpha_{k_{0}}-d_{0}\alpha_{k_{1}}}.
\een
\\

Type III:
$$ \Gamma_{ij}=\xy
(0,0); (10,0), **@{-};
(0,0)*+{\circ};(10,0)*+{\bullet};(5,2)*+{d};
(0,3)*+{i};(10,3)*+{j};
\endxy,
$$
where $1\leq i,j\leq l-1$, $i\neq j$. The contribution of $\Gamma_{ij}$ is
\ben
&&\frac{1}{d}\int_{\Mbar_{1,1}}\frac{\prod_{k=1,\neq i}^{l+1}\Lambda_{1}^{\vee}(\alpha_{i}-\alpha_{k})\cdot
\prod_{k=1}^{l-1}\Lambda_{1}^{\vee}(\alpha_{k}-\alpha_{i})\cdot\Lambda_{1}^{\vee}(\alpha_{l}+\alpha_{l+1}-2\alpha_{i}) }{\frac{\alpha_{i}-\alpha_{j}}{d}-\psi}\\
&&\cdot\frac{\frac{\alpha_{j}-\alpha_{i}}{d}\prod_{k=1}^{l-1}\prod_{a=1}^{d-1}(\alpha_{k}-\alpha_{j}+a\frac{\alpha_{j}-\alpha_{i}}{d})
\cdot\prod_{a=1}^{2d-1}(\alpha_{l}+\alpha_{l+1}-2\alpha_{j}+a\frac{\alpha_{j}-\alpha_{i}}{d})}{
(\frac{d!}{d^{d}})^{2}(\alpha_{i}-\alpha_{j})^d(\alpha_{j}-\alpha_{i})^d\prod_{k=1,\neq i,j}^{l+1}\prod_{a=0}^{d}(\alpha_{i}-\alpha_{k}+a\frac{\alpha_{j}-\alpha_{i}}{d})}\\
&=&\frac{(-1)^{l-1+(l-1)(d-1)}}{24d}
\prod_{k=1,\neq i,j}^{l-1}(\alpha_{i}-\alpha_{k})\cdot(\alpha_{l}-\alpha_{i})(\alpha_{l+1}-\alpha_{i})(\alpha_{l}+\alpha_{l+1}-2\alpha_{i})\\
&&\cdot\frac{
\prod_{a=1}^{2d-1}(\alpha_{l}+\alpha_{l+1}-2\alpha_{j}+a\frac{\alpha_{j}-\alpha_{i}}{d})}{
\prod_{k=1,\neq i,j}^{l-1}(\alpha_{j}-\alpha_{k})
\prod_{a=0}^{d}(\alpha_{i}-\alpha_{l}+a\frac{\alpha_{j}-\alpha_{i}}{d})(\alpha_{i}-\alpha_{l+1}+a\frac{\alpha_{j}-\alpha_{i}}{d})}.
\een
Setting $\alpha_{l}=\alpha_{l+1}=\alpha$, the above contribution is
\ben
&&\frac{(-1)^{(l-1)d}}{12d}
\prod_{k=1,\neq i,j}^{l-1}(\alpha_{i}-\alpha_{k})\cdot(\alpha-\alpha_{i})^{3}\\
&&\cdot\frac{
\prod_{a=1}^{2d-1}(2\alpha-2\alpha_{j}+a\frac{\alpha_{j}-\alpha_{i}}{d})}{
\prod_{k=1,\neq i,j}^{l-1}(\alpha_{j}-\alpha_{k})
\prod_{a=0}^{d}(\alpha_{i}-\alpha+a\frac{\alpha_{j}-\alpha_{i}}{d})^{2}}.
\een

Combining the three type of contributions, we obtain
\ben
N_{1,d}^{X}&=&\sum_{k_{0}=1}^{l-1}\frac{(-1)^{(l-1)d}(2d)!}{24(d!)^{2}d}\frac{\prod_{j=1,\neq k_{0}}^{l-1}(\alpha_{k_{0}}-\alpha_{j})}
{\prod_{j=1,\neq k_{0}}^{l-1}(\alpha-\alpha_{j})}
\Big(2d\sum_{a=1}^{d}\frac{1}{d+a}+1+
(\alpha-\alpha_{k_{0}})\sum_{j=1,\neq k_{0}}^{l-1}\frac{1}{\alpha-\alpha_{j}}\Big)\\
&&+\sum_{d_{0}+d_{1}=d}\sum_{k_{0}=1}^{l-1}\sum_{k_{1}=1}^{l-1}
\frac{(-1)^{(l-1)d+1}d_{0}(2d_{0})!(2d_{1})!}{24(d_{0}!)^{2}(d_{1}!)^{2}d_{1}}
\prod_{j=1,\neq k_{0}}^{l-1}(\alpha_{k_{0}}-\alpha_{j})\\
&&\cdot
\frac{1}
{\prod_{j=1,\neq k_{1}}^{l-1}(\alpha_{k_{1}}-\alpha_{j})}\cdot
\frac{\alpha-\alpha_{k_{0}}}{d\alpha-d_{1}\alpha_{k_{0}}-d_{0}\alpha_{k_{1}}}\\
&&+\sum_{i=1}^{l-1}\sum_{j=1,\neq i}^{l-1}\frac{(-1)^{(l-1)d}}{12d}
\prod_{k=1,\neq i,j}^{l-1}(\alpha_{i}-\alpha_{k})\cdot\frac{(\alpha-\alpha_{i})^{3}
\prod_{a=1}^{2d-1}(2\alpha-2\alpha_{j}+a\frac{\alpha_{j}-\alpha_{i}}{d})}{
\prod_{k=1,\neq i,j}^{l-1}(\alpha_{j}-\alpha_{k})
\prod_{a=0}^{d}(\alpha_{i}-\alpha+a\frac{\alpha_{j}-\alpha_{i}}{d})^{2}}\\
&=&\sum_{k_{0}=1}^{l-1}\frac{(-1)^{(l-1)d}(2d)!}{24(d!)^{2}d}\frac{\prod_{j=1,\neq k_{0}}^{l-1}(\alpha_{k_{0}}-\alpha_{j})}
{\prod_{j=1,\neq k_{0}}^{l-1}(\alpha-\alpha_{j})}
\Big(2d\sum_{a=1}^{d}\frac{1}{d+a}+1\Big)\\
&&+\sum_{i=1}^{l-1}\sum_{j=1,\neq i}^{l-1}\Big(\frac{(-1)^{(l-1)d}(2d)!}{24(d!)^{2}d}\frac{(\alpha_{i}-\alpha_{j})(\alpha-\alpha_{i})}{(\alpha-\alpha_{j})^{2}}\prod_{k=1,\neq i,j}^{l-1}\frac
{\alpha_{i}-\alpha_{k}}{\alpha-\alpha_{k}}
\Big)\\
&&+\sum_{k=1}^{l-1}\sum_{d_{0}+d_{1}=d}
\frac{(-1)^{(l-1)d+1}d_{0}(2d_{0})!(2d_{1})!}{24(d_{0}!)^{2}(d_{1}!)^{2}d_{1}d}\\
&&+\sum_{i=1}^{l-1}\sum_{j=1, \neq i}^{l-1}\sum_{d_{0}+d_{1}=d}
\frac{(-1)^{(l-1)d}d_{0}(2d_{0})!(2d_{1})!}{24(d_{0}!)^{2}(d_{1}!)^{2}d_{1}}
\prod_{k=1,\neq i,j}^{l-1}\frac{\alpha_{i}-\alpha_{k}}{\alpha_{j}-\alpha_{k}}\cdot
\frac{\alpha-\alpha_{i}}{d\alpha-d_{1}\alpha_{i}-d_{0}\alpha_{j}}\\
&&+\sum_{i=1}^{l-1}\sum_{j=1,\neq i}^{l-1}\frac{(-1)^{(l-1)d}}{12d}
\prod_{k=1,\neq i,j}^{l-1}\frac{\alpha_{i}-\alpha_{k}}{\alpha_{j}-\alpha_{k}}\cdot\frac{(\alpha-\alpha_{i})
\prod_{a=1}^{2d-1}(2\alpha-2\alpha_{j}+a\frac{\alpha_{j}-\alpha_{i}}{d})}{(\alpha-\alpha_{j})^{2}
\prod_{a=1}^{d-1}(\alpha_{i}-\alpha+a\frac{\alpha_{j}-\alpha_{i}}{d})^{2}}.
\een
Let us first assume $l\geq 3$. Note that $N_{1,d}^{X}$ is \emph{a priori} a rational number. So it is straightforward to see that, for fixed $1\leq i\neq j\leq l-1$, to cancel the denominators of the form $d\alpha-d_{1}\alpha_{i}-d_{0}\alpha_{j}$, it forces that there exist $b\in \mathbb{Q}$, $\beta_{ij}\in \mathbb{Q}\alpha_{1}+\cdots+\mathbb{Q}\alpha_{l-1}$, such that
\ben
b\alpha+\beta_{ij}&=&\sum_{d_{0}+d_{1}=d}
\frac{(-1)^{(l-1)d}d_{0}(2d_{0})!(2d_{1})!}{24(d_{0}!)^{2}(d_{1}!)^{2}d_{1}}
\cdot
\frac{(\alpha-\alpha_{j})^{2}}{d\alpha-d_{1}\alpha_{i}-d_{0}\alpha_{j}}\\
&&+\frac{(-1)^{(l-1)d}}{12d}
\cdot\frac{
\prod_{a=1}^{2d-1}(2\alpha-2\alpha_{j}+a\frac{\alpha_{j}-\alpha_{i}}{d})}{
\prod_{a=1}^{d-1}(\alpha_{i}-\alpha+a\frac{\alpha_{j}-\alpha_{i}}{d})^{2}}.
\een
Dividing both side by $\alpha$ and let $\alpha\rightarrow\infty$, we obtain
\ben
b&=&\sum_{d_{0}+d_{1}=d}
\frac{(-1)^{(l-1)d}d_{0}(2d_{0})!(2d_{1})!}{24(d_{0}!)^{2}(d_{1}!)^{2}d_{1}d}
+\frac{(-1)^{(l-1)d}2^{2d-1}}{12d}\\
&=&\frac{(-1)^{(l-1)d}}{24d}\Big(\sum_{d_{0}+d_{1}=d}
\frac{d_{0}(2d_{0})!(2d_{1})!}{(d_{0}!)^{2}(d_{1}!)^{2}d_{1}}
+4^{d}\Big).
\een
Then since
\ben
&&\sum_{j=1}^{l-1}\sum_{i=1,\neq j}^{l-1}\sum_{k=1,\neq i,j}^{l-1}\frac{\alpha_{i}-\alpha_{k}}{\alpha_{j}-\alpha_{k}}\\
&=&\sum_{i=1}^{l-1}\Big(\sum_{j=1,\neq i}^{l-1}(\alpha_{i}-\alpha_{j})\cdot
\sum_{j=1,\neq i}^{l-1}\frac{1}{-(\alpha_{j}-\alpha_{i})\prod_{k=1,\neq i,j}^{l-1}(\alpha_{j}-\alpha_{k})}\Big)\\
&=&\sum_{i=1}^{l-1}\Big(\sum_{j=1,\neq i}^{l-1}(\alpha_{i}-\alpha_{j})\cdot\frac{1}{\sum_{j=1,\neq i}^{l-1}(\alpha_{i}-\alpha_{j})}\Big)\\
&=&l-1,
\een
we have
\bea\label{70}
N_{1,d}^{X}&=&(l-1)b+(l-1)\sum_{d_{0}+d_{1}=d}
\frac{(-1)^{(l-1)d+1}d_{0}(2d_{0})!(2d_{1})!}{24(d_{0}!)^{2}(d_{1}!)^{2}d_{1}d}\nn\\
&=&\frac{(-1)^{(l-1)d}(l-1)4^{d}}{24d}.
\eea
For $l=2$, (\ref{70}) still holds, and has been proved in \cite{KP} without giving the details. Here we give another proof for this, which is interesting itself since we make use of the proof of the $l=3$ case to prove a combinatorial identity\footnote{Thanks Si-Qi Liu for telling the author that (\ref{72}) can also be \emph{proved} using Mathematica. }. It suffices to prove the following lemma.

\begin{lemma}\label{71}
\bea\label{72}
\frac{(2d)!}{(d!)^{2}}
\Big(2d\sum_{a=1}^{d}\frac{1}{d+a}+1\Big)
-\sum_{\substack{d_{0},d_{1}\geq 1\\d_{0}+d_{1}=d}}
\frac{d_{0}(2d_{0})!(2d_{1})!}{(d_{0}!)^{2}(d_{1}!)^{2}d_{1}}=4^{d}.
\eea
\end{lemma}
\emph{Proof}: Consider the case $l=3$. We have
\ben
N_{1,d}^{X}&=&-\frac{(2d)!}{24(d!)^{2}d}
\Big(2d\sum_{a=1}^{d}\frac{1}{d+a}+1\Big)\frac{(\alpha_{1}-\alpha_{2})^{2}}{(\alpha-\alpha_{1})(\alpha-\alpha_{2})}\\
&&+\frac{(2d)!}{24(d!)^{2}d}\Big(\frac{(\alpha_{1}-\alpha_{2})(\alpha-\alpha_{1})}{(\alpha-\alpha_{2})^{2}}
+\frac{(\alpha_{2}-\alpha_{1})(\alpha-\alpha_{2})}{(\alpha-\alpha_{1})^{2}}
\Big)\\
&&-2\sum_{d_{0}+d_{1}=d}
\frac{d_{0}(2d_{0})!(2d_{1})!}{24(d_{0}!)^{2}(d_{1}!)^{2}d_{1}d}\\
&&+\sum_{d_{0}+d_{1}=d}
\frac{d_{0}(2d_{0})!(2d_{1})!}{24(d_{0}!)^{2}(d_{1}!)^{2}d_{1}}
\Big(
\frac{\alpha-\alpha_{1}}{d\alpha-d_{1}\alpha_{1}-d_{0}\alpha_{2}}+\frac{\alpha-\alpha_{2}}{d\alpha-d_{1}\alpha_{2}-d_{0}\alpha_{1}}\Big)\\
&&+\sum_{i=1}^{2}\sum_{j=1,\neq i}^{2}\frac{(-1)^{(l-1)d}}{12d}
\cdot\frac{(\alpha-\alpha_{i})
\prod_{a=1}^{2d-1}(2\alpha-2\alpha_{j}+a\frac{\alpha_{j}-\alpha_{i}}{d})}{(\alpha-\alpha_{j})^{2}
\prod_{a=1}^{d-1}(\alpha_{i}-\alpha+a\frac{\alpha_{j}-\alpha_{i}}{d})^{2}}\\
&=&-\frac{(2d)!}{24(d!)^{2}d}
\Big(2d\sum_{a=1}^{d}\frac{1}{d+a}+1\Big)\frac{(\alpha_{1}-\alpha_{2})^{2}}{(\alpha-\alpha_{1})(\alpha-\alpha_{2})}
-2\sum_{d_{0}+d_{1}=d}
\frac{d_{0}(2d_{0})!(2d_{1})!}{24(d_{0}!)^{2}(d_{1}!)^{2}d_{1}d}\\
&&+\frac{1}{(\alpha-\alpha_{1})^{2}(\alpha-\alpha_{2})^{2}}\Big(
(\alpha-\alpha_{1})^{3}(b\alpha+\beta_{12}+\frac{(2d)!}{24(d!)^{2}d}(\alpha_{1}-\alpha_{2}))\\
&&+(\alpha-\alpha_{2})^{3}(b\alpha+\beta_{21}+\frac{(2d)!}{24(d!)^{2}d}(\alpha_{2}-\alpha_{1}))\Big).
\een
It forces that $\alpha-\alpha_{2}$ divides $b\alpha+\beta_{12}+\frac{(2d)!}{24(d!)^{2}d}(\alpha_{1}-\alpha_{2})$, and also
$\alpha-\alpha_{1}$ divides $b\alpha+\beta_{21}+\frac{(2d)!}{24(d!)^{2}d}(\alpha_{2}-\alpha_{1})$. Thus
\ben
\beta_{12}&=&(\frac{(2d)!}{24(d!)^{2}d}-b)\alpha_{2}-\frac{(2d)!}{24(d!)^{2}d}\alpha_{1},\\
\beta_{21}&=&(\frac{(2d)!}{24(d!)^{2}d}-b)\alpha_{1}-\frac{(2d)!}{24(d!)^{2}d}\alpha_{2},
\een
and
\ben
N_{1,d}^{X}&=&-\frac{(2d)!}{24(d!)^{2}d}
\Big(2d\sum_{a=1}^{d}\frac{1}{d+a}+1\Big)\frac{(\alpha_{1}-\alpha_{2})^{2}}{(\alpha-\alpha_{1})(\alpha-\alpha_{2})}
-2\sum_{d_{0}+d_{1}=d}
\frac{d_{0}(2d_{0})!(2d_{1})!}{24(d_{0}!)^{2}(d_{1}!)^{2}d_{1}d}\\
&&+\frac{b(\alpha-\alpha_{1})^{3}(\alpha-\alpha_{2})+b(\alpha-\alpha_{2})^{3}(\alpha-\alpha_{1})}{(\alpha-\alpha_{1})^{2}(\alpha-\alpha_{2})^{2}}\\
&=&-\frac{(2d)!}{24(d!)^{2}d}
\Big(2d\sum_{a=1}^{d}\frac{1}{d+a}+1\Big)\frac{(\alpha_{1}-\alpha_{2})^{2}}{(\alpha-\alpha_{1})(\alpha-\alpha_{2})}
-2\sum_{d_{0}+d_{1}=d}
\frac{d_{0}(2d_{0})!(2d_{1})!}{24(d_{0}!)^{2}(d_{1}!)^{2}d_{1}d}\\
&&+2b+\frac{b(\alpha_{1}-\alpha_{2})^{2}}{(\alpha-\alpha_{1})(\alpha-\alpha_{2})}.
\een
Therefore it forces that
\ben
b=\frac{(2d)!}{24(d!)^{2}d}
\Big(2d\sum_{a=1}^{d}\frac{1}{d+a}+1\Big).
\een
\pqed

\section{Integrality of $n_{1,d}$ for local Calabi-Yau 5-folds}
The Gopokumar-Vafa invariants $n_{0,d}(\gamma_{1},\cdots,\gamma_{k})$ for a Calabi-Yau $n$-fold $X$, where $\gamma_{1},\cdots,\gamma_{k}\in H^{*}(X)$ are defined by (see, e.g., \cite{KP}, \cite{PZ})
\bea\label{80}
\sum_{\beta\neq 0}\langle\gamma_{1},\cdots,\gamma_{k}\rangle_{0,k,\beta}^{X}q^{\beta}
=\sum_{\beta\neq 0}n_{0,\beta}(\gamma_{1},\cdots,\gamma_{k})\sum_{d=1}^{\infty}\frac{1}{d^{3-k}}q^{d\beta}.
\eea
When $n\geq 6$, the definition of Gopokumar-Vafa invariants in genus one\footnote{When $n\geq 4$, the Gromov-Witten invariants in genus at least two are triviall, due to the dimension constraint and the string equation.} is still absent. For $n=4$, the invariants $n_{1,d}$ are defined in \cite{KP}, and for $n=5$ in \cite{PZ}. The integrality of $n_{1,d}$ has been verified in low degrees in \cite{KP} for $X$ of the form (\ref{56}) when $n=4$ , and in \cite{PZ} the case $X=\mathrm{Tot}\big(\mathcal{O}(-1)^{\oplus 3}\rightarrow \mathbb{P}^{2}\big)$ when $n=5$. The remaining three cases for $n=5$ are $\mathcal{O}(-1)\oplus\mathcal{O}(-3)\rightarrow \mathbb{P}^{3}$,  $\mathcal{O}(-2)\oplus\mathcal{O}(-2)\rightarrow \mathbb{P}^{3}$,  $\mathcal{O}(-5)\rightarrow \mathbb{P}^{4}$. \\
For Calabi-Yau 5-folds, once we have $N_{1,d}$, $n_{0,i}(\gamma_{1})$ and $n_{0,i}(\gamma_{2},\gamma_{3})$ for $1\leq i\leq d$, all $\gamma_{1}\in H^{6}(X)$ and  all $\gamma_{2}, \gamma_{3}\in H^{4}(X)$ as inputs\footnote{We need also the Poincar\'{e} pairing on $H^{4}(X)\oplus H^{6}(X)$, which in the local cases are defined via the general principle mentioned in the footnote in Page 3. For example, for $X=K_{\mathbb{P}^{4}}$, we have $\langle H^{2},H^{3}\rangle^{X}=-\frac{1}{5}$.},  the invariants $n_{1,d}$ are defined through a complicated simultaneous recursion of many invariants. For the details we refer the reader to \cite{PZ}. The invariants $n_{0,i}(\gamma_{1})$ and $n_{0,i}(\gamma_{2},\gamma_{3})$ are defined by (\ref{80}), and the one-point and two-point genus zero Gromov-Witten invariants on the left of (\ref{80}) can be extracted from the formulae in \cite{Popa2} (see also \cite{GT}). Assuming the validity of our conjectural formulae (\ref{58}) and (\ref{59})  for $n=5$, we have checked the integrality of $n_{1,d}$ in for $1\leq d\leq 100$ for these three cases using a Maple programme, and for $1\leq d\leq 20$ we list them in the following.\\
\subsection{$X=K_{\mathbb{P}^{4}}$}
\ben
\langle H^{3}\rangle_{0,1,d}^{K_{\mathbb{P}^{4}}}=-\frac{1}{5}[x^{2}Q^{d}]\Big(e^{-xf(q)}\sum_{d\geq 0}q^{d}\frac{\prod_{s=0}^{5d-1}(-5x-s)}{\prod_{s=1}^{d}(x+s)^{5}}\Big),
\een
where $Q=qe^{f(q)}$ and the mirror map
\ben
f(q)=\sum_{d=1}^{\infty}q^{d}\frac{(-1)^{d}(5d)!}{d(d!)^{5}}.
\een

For $\langle H^{2},H^{2}\rangle_{0,2,d}^{K_{\mathbb{P}^{4}}}$, we follow the notations in the remark 3.4 in \cite{Popa2} and define $F(w,q)$ and $F_{i}(q)$ by
\ben
F(w,q)=\sum_{d=0}^{\infty}q^{d}\frac{\prod_{r=1}^{5d}(-5w-r)}{\prod_{r=1}^{d}(w+r)^{5}}=F(0,q)+\sum_{i=1}^{\infty}F_{i}(q)w^{i},
\een
and let
\ben
I_{1}(q)=1+q\frac{d}{dq}\frac{F_{1}(q)}{F(0,q)}.
\een
Then
\ben
\langle H^{2},H^{2}\rangle_{0,2,d}^{K_{\mathbb{P}^{4}}}=
-\frac{1}{5}[Q^{d}]\Bigg(-f(q)+\frac{\frac{F_{1}(q)}{F(0,q)}+q\frac{d}{dq}\frac{F_{2}(q)}{F(0,q)}}{I_{1}(q)}\Bigg).
\een

The conjectural formula (\ref{58}) in this case reads
\ben
\sum_{d=1}^{\infty}N_{1,d}Q^{d}=\frac{3}{8}f(q)-\frac{1}{8}\ln(1+5^{5}q)
-2\ln I_{1,1}(q)-\frac{1}{2}\ln I_{2,2}(q),
\een
where
\ben
I_{1,1}(q)&=&1+\sum_{d=1}^{\infty}\frac{(-1)^{d}(5d)!}{(d!)^{5}}q^{d},
\een
and
\ben
I_{2,2}(q)=1+\frac{1}{I_{1,1}(q)}\sum_{d=1}^{\infty}(-1)^{d}\frac{(nd)!}{(d!)^{n}}q^{d}
+\frac{1}{I_{1,1}(q)}\sum_{d=1}^{\infty}\Bigg(\frac{(-1)^{nd}nd(nd)!}{(d!)^{n}}\sum_{s=d+1}^{nd-1}\frac{1}{s}\Bigg)q^{d}.
\een

\begin{tabular}{ll}
\hline
$d$ & $n_{0,d}(H^{3})$ of $K_{\mathbb{P}^{4}}$ \\
\hline
1 & 130 \\
2 & -58345 \\
3 & 55837430\\
4 & -73589158000\\
5 & 115854201969950\\
6 & -204342355412313875\\
7 & 390051191739787697630 \\
8 & -789136006642194095804000\\
9 & 1669447288789130694933224250\\
10& -3658893431261650527639975955175\\
11& 8252627129183279407802045607394310\\
12& -19061509587415681611663858317767574480\\
13& 44917147949588887714507718293780333670230\\
14& -107667316864820156273192312584585440698457095\\
15& 261915168370711178492182001044618321338813469450\\
16& -645393917552138476376093839553201039666790189529280\\
17& 1608445644370011689169576347270893464407225594867091080\\
18& -4049011495564074654404411325327805800339427963862185528005\\
19& 10284566695008271699589128589728350347114600022600731093548670\\
20& -26334305024448861033964360994375819798940753071425109074393898000\\
\hline & \\
\end{tabular}\\

\begin{tabular}{ll}
\hline
$d$ & $n_{0,d}(H^{2},H^{2})$  of $K_{\mathbb{P}^{4}}$ \\
\hline
1 & 245\\
2 & -289035\\
3 & 499858460\\
4 & -1013558891950\\
5 & 2242341515096750\\
6 & -5241918236140466300\\
7 & 12728510402344664504790\\
8 & -31777727076990402350118750\\
9 & 81033105451821118038400330625\\
10& -210108099622343226675476798422750\\
11& 552168409753042747215570849250035965\\
12& -1467310058144521736953946444230597767540\\
13& 3935635344488399018105033615876566792311135\\
14& -10640097680708071622726078463500726511377961970\\
15& 28962447849234885737946426072924327337201062739625\\
16& -79304558059583206285015753062827770036005404547746270\\
17& 218282627813734979394266996554573007230698186101327263940\\
18& -603583836217572434350857017266811467968037123021501795632035\\
19& 1675850405011579381470981479391632353311673301482974031124005645\\
20& -4670116536398709153329791030555933322776990940284471434964091242500\\
\hline &\\
\end{tabular}\\

\begin{tabular}{ll}
\hline
$d$ & $n_{1,d}$  of $K_{\mathbb{P}^{4}}$ \\
\hline
1 & 0\\
2 & 0\\
3 & -27735575\\
4 & 138263175125\\
5 & -502345733521805\\
6 & 1625730914586631100\\
7 & -4991836999897827628150\\
8 & 14920114958100504172550700\\
9 & -43938600906882061090032617300\\
10& 128301145689055798368066779831220\\
11& -372790080292682641205105927773314550\\
12& 1080077789712734643768778488776735871550\\
13& -3124338094975833754174103026588005609926750\\
14& 9030728078747106325699133878410275027177326800\\
15& -26095962866424999144571647881038683645220729316310\\
16& 75414521979936074538799364083902488385523907361074200\\
17& -218002874312626664461636749643780260738215567343722123400\\
18& 630457897353302410598120441590533983765325492429949897687300\\
19& -1824211864568158295061578855171426326597137430136994284773704950\\
20& 5281330075502542531439277333474238318975916338665113159271700486035\\
\hline &\\
\end{tabular}

It is interesting to note that they are all multiples of 5, and when $5\nmid d$, $n_{1,d}$ is a multiple of 25.
\subsection{$X=\mathrm{Tot}\big(\mathcal{O}(-1)\oplus\mathcal{O}(-3)\rightarrow\mathbb{P}^{3}\big)$}
\ben
\langle H^{3}\rangle_{0,1,d}^{X}=\frac{(d-1)!(3d-1)!}{(d!)^{4}}.
\een

\ben
\langle H^{2},H^{2}\rangle_{0,2,d}^{X}=[q^{d}]\Bigg(\frac{\sum_{d=1}^{\infty}q^{d}\frac{(3d)!}{(d!)^{3}}\sum_{r=d+1}^{3d}\frac{1}{r}}
{1+\sum_{d=1}^{\infty}q^{d}\frac{(3d)!}{(d!)^{3}}}\Bigg).
\een

\ben
\sum_{d=1}^{\infty}N_{1,d}q^{d}=-\frac{1}{8}\ln(1-27q)-\frac{1}{2}\ln\Big(1+\sum_{d=1}^{\infty}q^{d}\frac{(3d)!}{(d!)^{3}}\Big).
\een

\begin{tabular}{llll}
\multicolumn{4}{c}{Table 1 : Low degree genus 0 and genus 1 BPS numbers of $\mathrm{Tot}\big(\mathcal{O}(-1)\oplus\mathcal{O}(-3)\rightarrow\mathbb{P}^{3}$\big)}\\[5pt]
\hline
$d$ &$n_{0,d}(H^{3})$&$n_{0,d}(H^{2},H^{2})$ & $n_{1,d}$   \\
\hline
1 &2&5& 0 \\
2 &7&53& 0 \\
3 &62&888& 135\\
4 &720&16578& 4069\\
5 &10090&336968& 102497\\
6 &158809&7208592& 2529330\\
7 &2714782&159953128& 62485370\\
8 &49299360&3644804226& 1549538856\\
9 &937750740&84757873392& 38632050468\\
10&18503320115& 2002782861068& 968230418446\\
11&376107425518&47940402636848& 24386703246083\\
12&7835027188272&1159841269631844& 616987529756004\\
13&166623467599342&28312447677391792& 15673085566208659\\
14&3606416097808937&696398907175066480& 399583442014671692\\
15&79251821904257590&17241740125645491096& 10220554875333281200\\
16&1764772740099673920&429315366375232815762& 262188626394087701664\\
17&39757622487694555282&10743399666271987545848& 6743753349276509395348\\
18&904958567371990915302&270039166920941445186084& 173872012409851929166786\\
19&20788888672249855553518&6814313281153255310131216& 4492655791971935260396097\\
20&481526012065391894029200&172564210354543917847594608& 116315885319017767137751283\\
\hline & \\
\end{tabular}

\subsection{$X=\mathrm{Tot}\big(\mathcal{O}(-2)\oplus\mathcal{O}(-2)\rightarrow\mathbb{P}^{3}\big)$ }
\ben
\langle H^{3}\rangle_{0,1,d}^{X}=\frac{(2d-1)!(2d-1)!}{(d!)^{4}}.
\een

\ben
\langle H^{2},H^{2}\rangle_{0,2,d}^{X}=[q^{d}]\Bigg(\frac{\sum_{d=1}^{\infty}q^{d}\frac{((2d)!)^{2}}{(d!)^{4}}\sum_{r=d+1}^{2d}\frac{1}{r}}
{1+\sum_{d=1}^{\infty}q^{d}\frac{((2d)!)^{2}}{(d!)^{4}}}\Bigg).
\een

\ben
\sum_{d=1}^{\infty}N_{1,d}q^{d}=-\frac{1}{8}\ln(1-16q)-\frac{1}{2}\ln \Big(1+\sum_{d=1}^{\infty}q^{d}\frac{((2d)!)^{2}}{(d!)^{4}}\Big).
\een

\begin{tabular}{llll}
\multicolumn{4}{c}{Table 2 : Low degree genus 0 and genus 1 BPS numbers of $\mathrm{Tot}\big(\mathcal{O}(-2)\oplus\mathcal{O}(-2)\rightarrow\mathbb{P}^{3}$\big)}\\[5pt]
\hline
$d$ & $n_{0,d}(H^{3})$ & $n_{0,d}(H^{2},H^{2})$ & $n_{1,d}$  \\
\hline
1 & 1 & 2  &0 \\
2 & 2 & 12& 0\\
3 & 11 & 122& 20\\
4 & 76 & 1344& 411\\
5 & 635& 16182& 6228\\
6 & 5926 & 204508& 92696\\
7 & 60095     &2683410& 1372416\\
8 & 647000     &36160512& 20351408\\
9 & 7296000   &497432288& 303008660\\
10& 85336790  &6954446148& 4529630140\\
11& 1028170055 &98509313850& 67986636924\\
12& 12695240996 &1410519352384& 1024271346252\\
13& 160018462071 &20380347529206& 15484823717804\\
14& 2052731611966 &296747545660524& 234834989626688\\
15& 26734938900985 &4349510282254174& 3571572918808416\\
16& 352829721754800 &64120438449094656& 54460621524782072\\
17& 4710828711092291 &950056145934862062& 832396434024038536\\
18& 63547901783133744 &14139866390015314240& 12750049354231063044\\
19& 865157668345976759 &211286868769225452618& 195680390778912132364\\
20& 11876040942305597380 &3168484757758896223680& 3008606422494946135414\\
\hline & & &\\
\end{tabular}

\hspace{1cm}\footnotesize{Department of Mathematical Sciences, Tsinghua University, Beijing, 100084, China }\\

\hspace{1cm}\footnotesize{\emph{E-mail address}: huxw08@mails.tsinghua.edu.cn}


\begin{thebibliography}{123}
\bibitem{ABK}Aganagic, Mina., Bouchard, Vincent., Klemm, Albrecht. \emph{Topological strings and (almost) modular forms.} Communications in Mathematical Physics 277.3 (2008): 771-819.
\bibitem{BCOV1}Bershadsky, M., Cecotti, S., Ooguri, H., Vafa, C. (1993). \emph{Holomorphic anomalies in topological field theories}. Nuclear Physics B, 405(2), 279-304.
\bibitem{BCOV2}Bershadsky, M., Cecotti, S., Ooguri, H.,  Vafa, C. (1994). \emph{Kodaira-Spencer theory of gravity and exact results for quantum string amplitudes}. Communications in Mathematical Physics, 165(2), 311-427.
\bibitem{GP}Graber, Tom., Pandharipande, Rahul. \emph{Localization of virtual classes}. Inventiones mathematicae 135.2 (1999): 487-518.
\bibitem{GT}Gholampour, Amin.,  Hsian-Hua Tseng. \emph{On computations of genus zero two-point descendant Gromov-Witten invariants}. arXiv preprint arXiv:1207.6071 (2012).
\bibitem{KP}Klemm, A., Pandharipande, R. \emph{Enumerative geometry of Calabi-Yau 4-folds.} Communications in Mathematical Physics 281.3 (2008): 621-653.
\bibitem{PZ}Pandharipande, Rahul.,  Zinger, Aleksey. \emph{Enumerative geometry of Calabi-Yau 5-folds.} arXiv preprint arXiv:0802.1640 (2008).
\bibitem{Popa1}Popa, Alexandra. \emph{The genus one Gromov-Witten invariants of Calabi-Yau complete intersections.} Transactions of the American Mathematical Society 365.3 (2013): 1149-1181.
\bibitem{Popa2}Popa, Alexandra. \emph{Two-point Gromov-Witten formulas for symplectic toric manifolds.} arXiv preprint arXiv:1206.2703 (2012).
\bibitem{Zinger1}Zinger, Aleksey. \emph{The reduced genus 1 Gromov-Witten invariants of Calabi-Yau hypersurfaces.} Journal of the American Mathematical Society 22.3 (2009): 691-737.
\end{thebibliography}
\end{document}